\documentclass[a4paper,11pt]{article}
\usepackage[margin=25mm,top=23mm,bottom=25mm]{geometry}
\usepackage{amsmath,amsthm}
\usepackage{fontspec,unicode-math}
\usepackage{xcolor,graphicx,tikz,caption}
\usepackage{placeins}
\usepackage[unicode,hidelinks]{hyperref}
\hypersetup{pdftitle={An Annular Spiral Construction for Anomalous Dissipation in the Forced Navier-Stokes Equations},
  pdfauthor={Changhong Li}}

\newtheorem{theorem}{Theorem}[section]
\newtheorem{proposition}[theorem]{Proposition}
\newtheorem{lemma}[theorem]{Lemma}
\numberwithin{equation}{section}
\newcommand{\T}{\mathbb T}
\allowdisplaybreaks[1]

\begin{document}
\begin{center}
{\LARGE\bfseries An Annular Spiral Construction for Anomalous Dissipation in the Forced Navier–Stokes Equations\par}
\vspace{.7em}
{\large Changhong Li\par}
\vspace{.2em}
{\small Email address: \href{mailto:lichanghong@amss.ac.cn}{lichanghong@amss.ac.cn}\par}
\end{center}

\section*{Abstract}

Kolmogorov's 1941 theory of turbulence (K41) predicts that the rate of kinetic energy dissipation in a turbulent flow can remain positive as viscosity tends to zero.
In 2023, Bruè and De Lellis used quasi-self-similar mixing to construct solutions of the forced three-dimensional Navier--Stokes equations exhibiting anomalous dissipation.
In this paper, we show that anomalous dissipation can still occur for a flow with spiral geometry in an annulus. In the horizontal plane, we wind a radial segment into a two-turn spiral and then apply alternating shears. After each pair of shears, the image of the segment shares a subarc of positive length with the original spiral. We compare the vertical velocity with sine profiles whose angular period halves with each pair. After the horizontal flow stops, diffusion damps the final profile. By embedding this planar construction into the $(2+\tfrac12)$-dimensional framework, we obtain smooth solutions of the three-dimensional equations with the same smooth initial datum. For a sequence of viscosities tending to zero, their dissipation over a fixed subannulus and one vertical period has a uniform positive lower bound, and their forces converge in $C^\alpha$ uniformly in time for every $0<\alpha<1$. We use this spiral construction to produce solutions of the forced three-dimensional Navier--Stokes equations that exhibit anomalous dissipation.

\section{Introduction}

\subsection{Background and motivation}

A central feature of three-dimensional turbulence is the energy cascade,
in which kinetic energy is transferred from large to small scales and is
ultimately dissipated by viscosity \cite{K41}.
We consider the forced incompressible Navier--Stokes equations with viscosity
$\mu_m>0$ on the three-dimensional torus $\T^3$,

\begin{equation}\label{eq:NS}
 \begin{cases}
 \partial_t u^m+u^m\cdot\nabla u^m+\nabla p^m
      =\mu_m\Delta u^m+f^m,\\
 \operatorname{div}u^m=0,\qquad u^m(0)=u_0,
 \end{cases}
\end{equation}
where $u^m$ is the velocity, $p^m$ the pressure, $f^m$ the external force,
and $u_0$ the initial velocity. For a sequence of viscosities $\mu_m\downarrow0$, the corresponding
solutions $u^m$ exhibit \emph{anomalous dissipation} on $[0,T]$ if the
kinetic energy dissipated by viscosity over this interval stays bounded
away from zero as $m\to\infty$, that is,

\begin{equation}\label{eq:ad-definition}
 \liminf_{m\to\infty}\mu_m\int_0^T\|\nabla u^m(t)\|_{L^2(\T^3)}^2\,dt>0 .
\end{equation}

DiPerna and Majda \cite[Sec.~2]{DM} described a
$(2+\tfrac12)$-dimensional construction for the Navier--Stokes equations.
Jeong and Yoneda \cite{JY21,JY} used this kind of flow to study
enhanced dissipation. Bru\`e and De Lellis \cite[Sec.~3]{BDL} used this framework to construct
solutions which satisfy the lower bound in \eqref{eq:ad-definition}
and thus exhibit anomalous dissipation. In their $(2+\frac12)$-dimensional construction,
$u^m=(v^m,\theta^m)$ and $f^m=(g^m,0)$ are independent of the third
coordinate. The vertical velocity solves the scalar equation

\begin{equation}\label{eq:scalar}
 \partial_t\theta^m+v^m\cdot\nabla\theta^m=\mu_m\Delta\theta^m,
 \qquad \theta^m(0)=\theta_0 .
\end{equation}

The viscous dissipation of $\theta^m$ contributes directly to
the dissipation term in \eqref{eq:ad-definition}. Bruè and De Lellis used the quasi-self-similar Peano snake construction of Alberti, Crippa and Mazzucato to prove anomalous dissipation for the forced three-dimensional Navier--Stokes equations \cite[Secs.~4 and~6]{BDL}. In this construction, the transported scalar takes both positive and negative values in successively smaller squares. Each step divides every square into
$5\times5$ smaller squares \cite[Sec.~8.1]{ACM}. At the end of the step, the transported scalar in each smaller
square is a rescaled copy of a profile supported on a straight or bent strip. If the initial square has side $L$, then after $n$ steps the transported
scalar integrates to zero over each of the resulting squares of side $L/5^n$
\cite[Secs.~6, 8.1 and~8.6]{ACM}.

In \cite{LI} we showed that anomalous dissipation persists for a class of small
normal perturbations of the generating curves satisfying the
area-preserving and equal-length conditions. In this paper, we construct a flow with spiral geometry in an annulus and prove that the corresponding solutions of the forced three-dimensional Navier--Stokes equations also exhibit anomalous dissipation. Our idea is to construct a spiral pattern like that of a tornado
viewed from above. After each pair of shears, the scalar profile used for comparison vanishes on the original spiral and its rotated copies. The points where these spirals meet a fixed ray are then half as far apart as before.

Our first idea was to wind a radial segment into a single spiral and
gradually decrease the distance between successive turns. Flohr and Vassilicos \cite[Sec.~2]{FV} studied the diffusion of a scalar patch wound around a vortex. In Appendix~\ref{sec:circular}, we construct a flow that winds a radial segment into a single spiral. Each particle stays on its initial circle, and the angular velocity is a time-dependent multiple of a fixed radial profile. On each circle this rotation shifts the scalar in angle and preserves its angular period. If the horizontal forces are uniformly bounded in $L^1(0,T;L^2(\T^2))$, the number of turns and the $L^2$ norm of the gradient of the scalar transported without diffusion stay bounded on $[0,T]$, uniformly in $m$. The energy estimate gives an upper bound $C_T\mu_m$ for the total viscous dissipation. This bound tends to zero as $\mu_m\to0$, so this construction cannot produce anomalous dissipation. To produce finer angular scales after the initial winding, we therefore allow particles to move between circles.

In the Peano construction the straight and bent strips join to form one strip \cite[Sec.~8.1]{ACM}. This led us to follow one material line, alternating a shear that changes the radius with an angular shear at fixed radius. The initial scalar changes sign across the line, whose image therefore lies in the zero set of the scalar transported without diffusion.

Alternating shears have been used to produce anomalous dissipation for passive scalars \cite{DEIJ,CCS} and for the forced Navier--Stokes equations \cite[Sec.~2]{BCCDS}. Elgindi and Liss used alternating shears of increasing spatial frequency to produce anomalous dissipation for a passive scalar. They noted that their construction also yields $(2+\tfrac12)$-dimensional Navier--Stokes solutions with anomalous dissipation \cite[Sec.~1]{EL}. Johansson and Sorella used the planar shears of Colombo, Crippa and Sorella to build a steady three-dimensional velocity field in their construction of anomalous dissipation for the forced four-dimensional Navier--Stokes equations \cite[Secs.~1 and~3.2]{JS}. Here the scalar is compared after each shear with a normalized spiral profile. The comparison controls the $L^2$ error at the time when the horizontal flow stops.

\subsection{Geometry and notation}\label{sec:notation}

Let $\T=\mathbb R/(2\pi\mathbb Z)$. We take integrals over $\T^j$ as
averages over one period in each coordinate. We write $\|\cdot\|_2$ and $\langle\cdot,\cdot\rangle$ for the spatial $L^2$ norm and inner product with the normalized measure above.

Write a point of $\T^3$ as
$(x',x_3)$, with $x'=(x_1,x_2)$. In the horizontal period cell
$(-\pi,\pi)^2$, we use polar coordinates $(r,\varphi)$ centered at
$(0,0)$, with $\varphi$ measured counterclockwise from the positive
$x_1$-axis.

Here $r$ is the distance to the origin. We construct the horizontal velocity in the annulus $\Omega=\{\sqrt{2}<r<2\}$ and estimate dissipation in its subannulus $\Omega'=\{3/2<r<\sqrt{15}/2\}$.

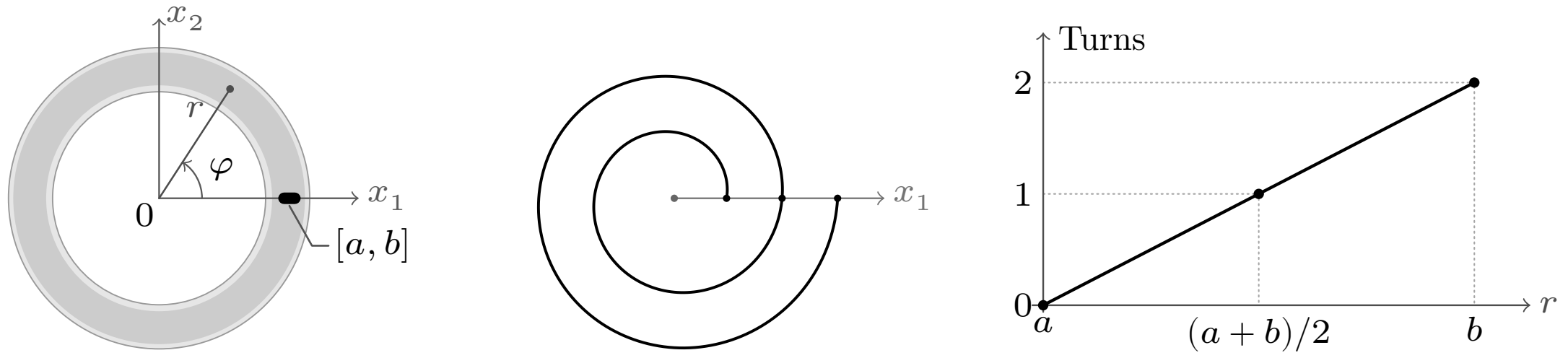
\begin{figure}[!htbp]
\centering
\begin{tikzpicture}[x=1cm,y=1cm,line cap=round,line join=round,font=\footnotesize]
  \pgfmathsetmacro{\segA}{sqrt(11)/2}
  \pgfmathsetmacro{\segB}{sqrt(13)/2}
  % Both plane views use equal horizontal and vertical units.
  \begin{scope}[x=.69cm,y=.69cm]
    \fill[even odd rule,black!10] (0,0) circle[radius=2] (0,0) circle[radius={sqrt(2)}];
    \fill[even odd rule,black!20] (0,0) circle[radius={sqrt(15)/2}] (0,0) circle[radius=1.5];
    \foreach \rad in {1.41421356,2}
      \draw[black!40,line width=.35pt] (0,0) circle[radius=\rad];
    \draw[->,black!65,line width=.44pt] (0,0)--(2.65,0) node[right,inner sep=2pt] {$x_1$};
    \draw[->,black!65,line width=.44pt] (0,0)--(0,2.39) node[right,inner sep=1.5pt] {$x_2$};
    \node[below left,inner sep=1pt] at (0,0) {$0$};
    \draw[black!70,line width=.48pt] (0,0)--(57:1.73)
      node[pos=.70,above left,inner sep=1pt] {$r$};
    \fill[black!70] (57:1.73) circle[radius=1pt];
    \draw[->,black!65,line width=.43pt] (.57,0) arc(0:57:.57);
    \node at (25:.91) {$\varphi$};
    \draw[black,line width=3pt] (\segA,0)--(\segB,0);
    \foreach \rad in {\segA,\segB} \fill[black] (\rad,0) circle[radius=1.03pt];
    \draw[black!68,line width=.5pt] (1.73,-.10)--(2.03,-.63)--(2.25,-.63);
    \node[anchor=west,inner sep=1.4pt,text=black] at (2.25,-.63) {$[a,b]$};
  \end{scope}
  \begin{scope}[shift={(4.72,0)},x=.83cm,y=.83cm]
    % The schematic widens the radial span only in this panel.
    \pgfmathsetmacro{\drawingInner}{.32*\segB}
    \draw[->,black!55,line width=.43pt] (0,0)--(2.32,0) node[right,inner sep=2pt] {$x_1$};
    \fill[black!60] (0,0) circle[radius=.95pt];
    \draw[black,line width=.74pt]
      plot[domain=0:720,samples=361,variable=\ang]
      ({(\drawingInner+(\segB-\drawingInner)*\ang/720)*cos(\ang)},
       {(\drawingInner+(\segB-\drawingInner)*\ang/720)*sin(\ang)});
    \foreach \s in {0,.5,1}
      \fill[black] ({\drawingInner+(\segB-\drawingInner)*\s},0) circle[radius=.95pt];
  \end{scope}
  \begin{scope}[shift={(8.10,-.98)},x=3.95cm,y=1.02cm]
    \foreach \xx/\hh in {.5/1,1/2} {
      \draw[black!35,densely dotted] (0,\hh)--(\xx,\hh)--(\xx,0);
    }
    \draw[->,black!70,line width=.45pt] (0,-.06)--(0,2.45);
    \node[anchor=west,inner sep=2pt] at (.015,2.40) {Turns};
    \draw[->,black!70,line width=.45pt] (-.025,0)--(1.13,0) node[right,inner sep=2pt] {$r$};
    \foreach \hh in {0,1,2}
      \node[anchor=east,inner sep=2.5pt] at (0,\hh) {$\hh$};
    \draw[black,line width=.82pt] (0,0)--(1,2);
    \foreach \xx/\hh/\lab in {0/0/a,.5/1/{(a+b)/2},1/2/b} {
      \fill[black] (\xx,\hh) circle[radius=1.35pt];
      \node[below,inner sep=2.8pt] at (\xx,0) {$\lab$};
    }
  \end{scope}
\end{tikzpicture}
\caption[The marked segment and the two-turn spiral.]{In the horizontal plane, the bold segment in the shaded subannulus $\Omega'$ winds into a two-turn spiral. The middle panel enlarges the spacing between turns for clarity. The graph shows the number of turns at radius $r$.}
\label{fig:annular-geometry}
\end{figure}

At time zero, we mark the particles on the segment of the positive $x_1$-axis from $(a,0)$ to $(b,0)$, with $a=\sqrt{11}/2$ and $b=\sqrt{13}/2$. We use this segment as the initial material line (Figure~\ref{fig:annular-geometry}). The initial scalar vanishes along it.

During the initial winding, each marked particle rotates about the origin at its initial radius $r$. At the end of this initial winding, the number of counterclockwise turns is $2(r-a)/(b-a)$, so it increases linearly with $r$. For example, the inner endpoint stays fixed, the midpoint makes one full turn, and the outer endpoint makes two. In this construction, all three points end on the positive $x_1$-axis. The image of the segment is then a two-turn spiral.

\FloatBarrier

\subsection{Main result}

In this subsection, we state the main theorem of this paper.

\begin{theorem}\label{thm:main}

Fix $0<T_*<1$. There exist a sequence
$\mu_m\downarrow0$ and smooth solutions $u^m=(v^m,\theta^m)$ of
\eqref{eq:NS} on $[0,T_*]\times\T^3$, with the same smooth initial datum
$u_0=(0,0,\theta_0)$ and smooth forces $f^m=(g^m,0)$. Both $u^m$ and $f^m$ are independent of $x_3$, and
$\int_{\T^2}\theta_0\,dx'=0$, $\|\theta_0\|_2=1$.

\begin{enumerate}

\renewcommand{\labelenumi}{(\arabic{enumi})}

\item

For each integer $m\ge0$, the horizontal velocity $v^m$ is supported in $\Omega$. Its flow first maps the marked segment onto the two-turn spiral of Section~\ref{sec:notation}. After each of the $m+1$ pairs of shears, the image of the segment shares a subarc of positive length with this spiral between radii $a$ and $b$.

\item

There exists a force $f$ such that, for every $0<\alpha<1$,

\[
 f^m\longrightarrow f\quad\text{in }C([0,T_*];C^\alpha(\T^3)).
\]

The solutions exhibit anomalous dissipation in $\Omega'\times\T$, with

\begin{equation}\label{eq:main-local}
 \liminf_{m\to\infty}\mu_m\int_0^{T_*}\!\int_{\Omega'\times\T}
       |\nabla u^m(s)|^2\,dx'\,dx_3\,ds\ge\frac38 .
\end{equation}

\end{enumerate}

\end{theorem}

\textbf{Organization.} Section~\ref{sec:geometry} constructs the horizontal flow and follows the marked segment. Section~\ref{sec:half} estimates the error from one shear. Section~\ref{sec:scalar-dissipation} sums these errors and proves the scalar dissipation bound. Section~\ref{sec:NS} proves convergence of the forces and completes the proof of Theorem~\ref{thm:main}. We finally introduce a single-spiral model with a fixed angular velocity profile in Appendix~\ref{sec:circular}.

\section{The annular construction}\label{sec:geometry}

Fix $T_*$ as in Theorem~\ref{thm:main}. We construct the horizontal velocity and follow the radial segment marked in Section~\ref{sec:notation}.

\begin{proposition}[Spiral geometry]\label{prop:geometry}
For the times and shear parameters in Sections~\ref{sec:preparation} and~\ref{sec:shears}, the construction gives smooth divergence-free velocities $v^m$ on $[0,T_*]\times\T^2$ with $v^m(0)=0$.
For each integer $m\ge0$, the horizontal velocity $v^m$ is supported in $\Omega$.
The initial winding maps the marked segment onto the two-turn spiral of Section~\ref{sec:notation}. After each of the $m+1$ pairs of shears, the image of the segment contains a subarc of positive length on this spiral between radii $a$ and $b$. The velocity vanishes after the last pair.

\end{proposition}

\subsection{Coordinates and the initial scalar}

We choose a radial coordinate $\sigma$ with $d\sigma=r\,dr$ to simplify area integrals. Let $y$ denote the polar angle $\varphi$ measured in turns, with one turn equal to $2\pi$ radians. The two coordinates are given by

\[
 \sigma=\frac{r^2}{2},\qquad y=\frac{\varphi}{2\pi}\in\mathbb R/\mathbb Z .
\]

Adding an integer to $y$ amounts to a whole number of turns and thus gives the same point on the circle. For example, in this paper both $y=0$ and $y=1$ describe the positive $x_1$-axis.

Let $J$ denote the interval of $\sigma$ occupied by the marked segment as $r$ varies from $a$ to $b$. In this coordinate the two annuli and the marked interval are given by

\[
 \Omega=\{1<\sigma<2\},\qquad \Omega'=\{9/8<\sigma<15/8\},\qquad J=[11/8,13/8].
\]

Measured in $\sigma$, the distances from $J$ to the inner and outer boundaries of $\Omega'$ are $11/8-9/8=1/4$ and $15/8-13/8=1/4$, respectively.

The period cell $(-\pi,\pi)^2$ has Euclidean area $(2\pi)^2$. We divide Euclidean area by $(2\pi)^2$ so that the whole cell has normalized area one. As the outer radius of $\Omega$ is $2<\pi$, the annulus lies entirely inside this cell. Using $d\sigma=r\,dr$ and $d\varphi=2\pi\,dy$, we obtain

\[
 dx'=\frac{dx_1\,dx_2}{(2\pi)^2}
    =\frac{r\,dr\,d\varphi}{(2\pi)^2}
    =\frac{1}{2\pi}\,d\sigma\,dy \qquad \text{on }\Omega .
\]

Along the two-turn spiral in Figure~\ref{fig:annular-geometry}, the angle increases by two turns as the radius increases from $a$ to $b$. Let $\lambda$ denote this increase in turns per unit radius. At radius $r$ the angle of the spiral is $\lambda(r-a)$ turns. We use $\Psi$ for the angular coordinate measured from this angle, and use the same formula throughout $\Omega$. Set

\[
 \lambda=\frac{2}{b-a},\qquad \Psi=y-\lambda(r-a)\in\mathbb R/\mathbb Z .
\]

For $a\le r\le b$, the original spiral consists of the points with $\Psi=0$. Fix an integer $q\ge1$. Starting with the original spiral, rotate it about the origin by $1/q$ of a turn at each of $q-1$ steps. Each step adds $1/q$ to $y$ at fixed $r$, so on the $q$ spirals we have $\Psi=0,1/q,\ldots,(q-1)/q$, respectively. Their union is given by

\begin{equation}\label{eq:rotated-spirals}
 q\Psi=0\quad\text{in }\mathbb R/\mathbb Z,\qquad a\le r\le b .
\end{equation}

For example, when $q=2$, both $\Psi=0$ and $\Psi=1/2$ satisfy \eqref{eq:rotated-spirals}. At each radius, the point with $\Psi=1/2$ is half a turn from the point with $\Psi=0$ on the original spiral. As $r$ varies from $a$ to $b$, the points with $\Psi=1/2$ trace the spiral rotated by half a turn.

Follow a fixed ray outward and let $r_1<r_2$ be the radii of two consecutive intersections with the spirals. The angle $y$ stays fixed. By \eqref{eq:rotated-spirals}, there is an intersection each time $\lambda(r-a)-y$ reaches a multiple of $1/q$. This quantity increases linearly with $r$, so it increases by $1/q$ from $r_1$ to $r_2$. Thus

\[
 \lambda(r_2-r_1)=\frac1q,\qquad
 r_2-r_1=\frac{1}{q\lambda}=\frac{b-a}{2q}.
\]

When $q$ doubles, the radial spacing becomes $(b-a)/(4q)$. The new intersection is at radius $(r_1+r_2)/2$.

Let $A(s)$ denote the sine value for an angle of $s$ turns. For the initial scalar take $s=y$. On the positive $x_1$-axis $A(0)=0$, and a quarter turn counterclockwise gives $A(1/4)=1$. Fix a radius $r$. As the angle $y$ varies around the circle, the sine value $A(y)$ changes. The factor $\chi(r^2/2)$ multiplying $A(y)$ stays the same because it depends only on $r$. We choose $\chi\ge0$ everywhere and require $\chi(\sigma)>0$ for $\sigma\in J$, which is the marked interval. We also require $\chi$ to vanish near both the inner and outer boundaries of $\Omega$.

\[
 A(s)=\sin(2\pi s),\qquad \chi\in C_c^\infty((5/4,7/4)),\qquad \chi\ge0,\quad \min_J\chi>0 .
\]

Define the initial scalar $\theta_0$ on $\Omega$ by multiplying these factors,

\begin{equation}\label{eq:initial-profile}
 \theta_0(x')=\chi(\sigma)A(y).
\end{equation}

Extend $\theta_0$ to $\T^2$ by setting it equal to zero on $\T^2\setminus\Omega$. This extension is smooth because $\chi$ vanishes near both boundary circles. Since $\int_0^1 A(y)\,dy=0$, the spatial mean of $\theta_0$ is zero. It changes sign across the marked segment because $\chi>0$ on $J$ and $A(y)=\sin(2\pi y)$ changes sign at $y=0$.

We scale $\chi$ so that $\|\theta_0\|_2=1$. Using the normalized area element and $\int_0^1 A(y)^2\,dy=1/2$, we obtain

\[
 \|\theta_0\|_2^2
 =\frac{1}{2\pi}\int_1^2\chi(\sigma)^2\,d\sigma\int_0^1 A(y)^2\,dy
 =\frac{1}{4\pi}\int_1^2\chi(\sigma)^2\,d\sigma .
\]

Choosing $\int_1^2\chi(\sigma)^2\,d\sigma=4\pi$ thus gives unit $L^2$ norm and fixes the initial scalar energy $\tfrac12\|\theta_0\|_2^2$ at $\tfrac12$.

\subsection{Winding the radial segment}\label{sec:preparation}

Let $t_0$ be the end of the initial winding. For $n\ge0$, shear pair $n$ starts at $t_n$ and ends at $t_{n+1}$. Choose the times so that

\[
 0<t_0<t_1<\cdots,\qquad t_n\uparrow T_* .
\]

To make each motion smooth in time, choose a nonnegative function $\eta$ with integral one,

\[
 \eta\in C_c^\infty((0,1)),\qquad \eta\ge0,\qquad \int_0^1\eta(t)\,dt=1 .
\]

Extend $\eta$ to $\mathbb R$ by setting $\eta(t)=0$ for $t\notin(0,1)$.

On a time interval $I=(t_I,t_I+s)$ with starting time $t_I$ and duration $s>0$, let $\eta_I$ give the rate at which the prescribed motion is applied. Its integral $\tau_I(t)$ is the fraction completed by time $t$. Define

\[
 \eta_I(t)=s^{-1}\eta((t-t_I)/s),\qquad
 \tau_I(t)=\int_{t_I}^t\eta_I(\xi)\,d\xi .
\]

Thus $\tau_I(t_I)=0$ and $\tau_I(t_I+s)=1$. To localize the velocity in the annulus, choose a radial cutoff $\zeta$ equal to one on an interval extending $[5/4,7/4]$ by $1/16$ at each end,

\[
 \zeta\in C_c^\infty((9/8,15/8)),\qquad 0\le\zeta\le1,\qquad \zeta=1\quad (19/16\le\sigma\le29/16).
\]

A dot over a physical quantity denotes its time derivative along a particle trajectory. On the winding interval $I=(0,t_0)$, define the horizontal velocity $v$ by prescribing the radial and angular motion of each particle in $\Omega$ and setting $v=0$ outside $\Omega$. Then we have

\[
 \dot \sigma=0,\qquad \dot y=\eta_I(t)\zeta(\sigma)\lambda(r-a).
\]

The velocity is divergence free since $\partial_\sigma\dot\sigma+\partial_y\dot y=0$. On $\operatorname{supp}\chi$, each particle keeps its radius and turns through $\tau_I(t)\lambda(r-a)$ revolutions by time $t$. Let $\rho$ denote the scalar transported from $\theta_0$ by $v$ without diffusion. On $\operatorname{supp}\chi$, the initial angle is $y-\tau_I(t)\lambda(r-a)$. Thus

\[
 \rho(t,x')=\chi(\sigma)A\bigl(y-\tau_I(t)\lambda(r-a)\bigr)
 \qquad (x'\in\Omega).
\]

For $\mu>0$, $\theta$ is the solution of the advection--diffusion equation with the velocity $v$ constructed above and diffusivity $\mu$. Its initial value is the same as that of the transported scalar $\rho$,

\[
 \theta(0)=\rho(0)=\theta_0.
\]

At the end of the winding, $\tau_I(t_0)=1$. Substituting this in the formula for $\rho$ gives

\[
 \rho(t_0)=\chi(\sigma)A(\Psi).
\]

We use the method of Drivas, Elgindi, Iyer and Jeong \cite{DEIJ} to compare the viscous and transported scalars.

Consider any interval $[t_I,t_I+s]$. Let $\theta$ solve the advection--diffusion equation with diffusivity $\mu$ and let $\rho$ solve the transport equation with the same smooth divergence-free velocity $v$. We allow different initial values at $t_I$. Subtracting the equations gives

\[
 w=\theta-\rho,\qquad
 \partial_t w+v\cdot\nabla_{x'}w=\mu\Delta_{x'}\theta .
\]

Taking the $L^2(\T^2)$ inner product with $w$, using $\operatorname{div}_{x'}v=0$ and integrating by parts in $x'$, we obtain

\[
 \begin{aligned}
 \frac{d}{dt}\left(\frac{\|w\|_2^2}{2}\right)
 &=-\mu\langle\nabla_{x'}\theta,\nabla_{x'}(\theta-\rho)\rangle\\
 &=-\mu\|\nabla_{x'}\theta\|_2^2
      +\mu\langle\nabla_{x'}\theta,\nabla_{x'}\rho\rangle\\
 &\le -\frac{\mu}{2}\|\nabla_{x'}\theta\|_2^2
      +\frac{\mu}{2}\|\nabla_{x'}\rho\|_2^2 .
 \end{aligned}
\]

Integrating from $t_I$ to $t\in[t_I,t_I+s]$ yields

\begin{equation}\label{eq:prep}
 \begin{aligned}
 &\|\theta(t)-\rho(t)\|_2^2
       +\mu\int_{t_I}^t\|\nabla_{x'}\theta(\xi)\|_2^2\,d\xi\\
 &\qquad\le\|\theta(t_I)-\rho(t_I)\|_2^2
       +\mu\int_{t_I}^t\|\nabla_{x'}\rho(\xi)\|_2^2\,d\xi .
 \end{aligned}
\end{equation}

During the initial winding, $\theta(0)=\rho(0)$ and the explicit formula for $\rho$ is smooth and independent of $\mu$. Hence \eqref{eq:prep} gives

\[
 \|\theta(t_0)-\chi(\sigma)A(\Psi)\|_2
 \le\left(\mu\int_0^{t_0}\|\nabla_{x'}\rho(t)\|_2^2\,dt\right)^{1/2}
 \le C\sqrt\mu .
\]

\subsection{One pair of shears}\label{sec:shears}

After the initial winding the scalar transported without diffusion is $\chi(\sigma)A(\Psi)$. One pair of shears is designed to double the number of periods of the comparison profile on each circle where $\chi>0$. A rotation depending only on the radius preserves angular spacing on each circle. We therefore start with a shear that changes the radius.

At the start of pair $n$, let $q$ be the number of periods of the comparison profile around a circle where $\chi>0$. Let $p$ be the coefficient of $\sigma$ in the argument of $A$. The ratio $K=p/q$ controls the radial displacement. For Theorem~\ref{thm:main}, we choose one sufficiently large even integer $K\ge16$ after \eqref{eq:stoperror} and use it in every pair. We take

\[
 q=2^n,\qquad p=Kq,\qquad n=0,1,2,\ldots .
\]

To compute the changes of phase in the two shears separately, we introduce the angular coordinate $z$. At a fixed radius it measures the angle in turns from the point where $\Psi=K\sigma$ on the circle. Set

\[
 z=\Psi-K\sigma\in\mathbb R/\mathbb Z .
\]

At fixed $\sigma$, the coordinates $z$ and $y$ differ only by a constant, namely $z-y=-\lambda(\sqrt{2\sigma}-a)-K\sigma$. Since $z=y-\lambda(r-a)-K\sigma$ and $r=\sqrt{2\sigma}$, we have

\[
 \frac{dr}{d\sigma}=\frac1r,\qquad
 \frac{\partial(\sigma,z)}{\partial(\sigma,y)}
 =\det\begin{pmatrix}
     1 & 0\\[2pt] -\lambda/r-K & 1
 \end{pmatrix}
 =1.
\]

In particular,

\[
 \Psi=K\sigma+z,\qquad q\Psi=p\sigma+qz,\qquad
 dx'=\frac1{2\pi}\,d\sigma\,dz .
\]

For positive integers $p,q$, define the comparison function

\[
 \phi_{p,q}(x')=
 \begin{cases}
   \chi(\sigma)A(p\sigma+qz), & x'\in\Omega,\\
   0, & x'\in\T^2\setminus\Omega .
 \end{cases}
\]

At fixed $\sigma$, the angular period is $1/q$ turns wherever $\chi(\sigma)>0$. Since $\int_0^1 A(p\sigma+qz)^2\,dz=1/2$, the normalization of $\chi$ gives $\|\phi_{p,q}\|_2=1$.

To define a smooth periodic displacement function, we modify $s\mapsto s$ only near the endpoints of $[-1/2,1/2]$. Using $\eta$ from Section~\ref{sec:preparation}, define the cutoff

\[
 M_\gamma(s)=\int_{-\infty}^{(1/2-|s|)/\gamma-1}\eta(t)\,dt,
 \qquad |s|\le\tfrac12,\quad 0<\gamma<\tfrac18 .
\]

Then $M_\gamma=1$ for $|s|\le1/2-2\gamma$ and $M_\gamma=0$ for $1/2-\gamma\le|s|\le1/2$. Repeat the values of $M_\gamma$ on $[-1/2,1/2]$ every unit to define it on the whole line, so $M_\gamma(s+1)=M_\gamma(s)$. For pair $n$, take $h=h_n$. Keep the displacement equal to $s$ on $[-1/2+h,1/2-h]$. On the two intervals of width $h$ next to the endpoints, use $M_{h/2}$ to make the displacement vanish smoothly near each end. Set

\[
 S_h(s)=s\,M_{h/2}(s)\quad (|s|\le\tfrac12),\qquad
 S_h(s+1)=S_h(s),\qquad 0<h=h_n<1/64 .
\]

Since $M_{h/2}$ is even and vanishes near both endpoints, $S_h$ is smooth and odd, with $|S_h|\le1/2$. For every integer $k$,

\[
 \int_0^1 S_h(s)\,ds=0,\qquad
 S_h(s)=s-k,\quad S_h^{\prime}(s)=1\qquad (|s-k|\le\tfrac12-h).
\]

Differentiating $M_\gamma$ and applying the product rule to $sM_{h/2}(s)$ give, for $j\ge1$,

\[
 \|M_\gamma^{(j)}\|_\infty\le\gamma^{-j}\|\eta^{(j-1)}\|_\infty,
 \qquad \|S_h^{(j)}\|_\infty\le C_jh^{-j},
\]

where $C_j$ depends only on $j$ and $\eta$. Let $R$ and $T$ denote the first and second shear maps. For particles that stay in $19/16\le\sigma\le29/16$, set

\begin{equation}\label{eq:shears}
 R(\sigma,z)=\left(\sigma-\frac{S_h(qz)}p,z\right),\qquad
 T(\sigma,z)=\left(\sigma,z-\frac{S_h(p\sigma)}{2q}\right).
\end{equation}

In \eqref{eq:shears}, we keep $z$ fixed during the radial shear $R$, so changing the radius also changes the angle $y=z+K\sigma+\lambda(r-a)$. During the angular shear $T$, we keep the radius fixed.

To keep the first shear divergence free where the radial cutoff $\zeta$ varies, define

\[
 Q_h(s)=\int_0^s S_h(\xi)\,d\xi .
\]

Since $S_h$ integrates to zero over one period, $Q_h$ is periodic. On the time interval $I$ of either shear, set

\begin{equation}\label{eq:ham}
 (\dot\sigma,\dot z)=\eta_I(t)
 \begin{cases}
 \bigl(-\zeta(\sigma)S_h(qz)/p,\ \zeta'(\sigma)Q_h(qz)/(pq)\bigr) & (R),\\[2pt]
 \bigl(0,-\zeta(\sigma)S_h(p\sigma)/(2q)\bigr) & (T).
 \end{cases}
\end{equation}

Since $dx'=(2\pi)^{-1}d\sigma\,dz$ and $Q_h'=S_h$, the divergence is

\[
 \operatorname{div}_{x'}v=\partial_\sigma\dot\sigma+\partial_z\dot z
 =\begin{cases}
 \eta_I(t)\zeta'(\sigma)\bigl(-S_h(qz)+Q_h'(qz)\bigr)/p=0 &(R),\\
 0 &(T).
 \end{cases}
\]

On $[19/16,29/16]$, $\zeta=1$ and $\zeta'=0$. Integrating \eqref{eq:ham} over $I$ with $\int_I\eta_I=1$ gives \eqref{eq:shears}.

Let $e_r$ and $e_\varphi$ be the outward radial and counterclockwise tangential unit vectors. In the horizontal plane, $\dot r=\dot\sigma/r$ and $\dot\varphi=2\pi(\dot z+K\dot\sigma+\lambda\dot r)$ give

\[
 v=\dot r\, e_r+r\dot\varphi\,e_\varphi
   =\frac{\dot\sigma}{r}e_r+
       2\pi r\bigl(\dot z+(\lambda/r+K)\dot\sigma\bigr)e_\varphi .
\]

The angular factors are periodic, and $\zeta$ vanishes near both boundary circles. Thus $v$ extends smoothly by zero outside $\Omega$. By \eqref{eq:ham} and $|S_h|\le1/2$, the change in $\sigma$ during pair $n$ is at most $1/(2K2^n)$. Through pair $m$, this gives

\begin{equation}\label{eq:radial-margin}
 |\sigma(t)-\sigma(t_0)|
 \le \sum_{n=0}^m\frac{1}{2K2^n}<\frac1K
 \qquad (t_0\le t\le t_{m+1}) .
\end{equation}

For $\sigma(t_0)\in\operatorname{supp}\chi\subset(5/4,7/4)$, the distance to either end of $[19/16,29/16]$ is at least $1/16$. By \eqref{eq:radial-margin}, $|\sigma(t)-\sigma(t_0)|<1/K\le1/16$, so $19/16<\sigma(t)<29/16$. Hence $\zeta(\sigma(t))=1$ throughout the shears.

Let $(\sigma,z)$ be the position of such a particle at the end of a shear. By \eqref{eq:shears}, we can obtain the starting point by adding $S_h(qz)/p$ to $\sigma$ in $R$, or $S_h(p\sigma)/(2q)$ to $z$ in $T$. Set $s=qz$ in $R$ and $s=p\sigma$ in $T$. If $|s-k|\le1/2-h$ for some integer $k$, then $S_h(s)=s-k$. Substituting the starting point in the sine factor and using the periodicity of $A$, we obtain

\[
 \begin{aligned}
 R:\quad A\!\left(p\left(\sigma+\frac{S_h(qz)}p\right)+qz\right)
 &= A(p\sigma+qz+S_h(qz))=A(p\sigma+2qz),\\
 T:\quad A\!\left(p\sigma+2q\left(z+\frac{S_h(p\sigma)}{2q}\right)\right)
 &= A(p\sigma+2qz+S_h(p\sigma))=A(2p\sigma+2qz).
 \end{aligned}
\]

At the end of pair $n$, the comparison function on $\Omega$ is $\phi_{2p,2q}=\chi(\sigma)A(2q\Psi)$. Its sine factor has an angular period of $1/(2q)$. Since $2q=2^{n+1}$ and $2p=K(2q)$, pair $n+1$ starts with this same function. In Section~\ref{sec:half} we estimate the change in $\chi$ under $R$ and the error near the ends of each period in both shears.

Fix finitely many shear pairs and their parameters. The transported scalar $\rho$ is then smooth and independent of $\mu$. Using $\theta(0)=\rho(0)$ in \eqref{eq:prep}, we obtain

\[
 \mu\int_0^{T_*}\|\nabla_{x'}\theta(t)\|_2^2\,dt
 \le \mu\int_0^{T_*}\|\nabla_{x'}\rho(t)\|_2^2\,dt
 \longrightarrow 0 \qquad (\mu\to0).
\]

We then let the number of pairs increase as $\mu_m\to0$ in Section~\ref{sec:scalar-dissipation}.

For a scalar function $\phi(\sigma,z)$, using $\partial_r\sigma=r$, $\partial_r z=-\lambda-Kr$ and $\partial_\varphi z=1/(2\pi)$, we obtain

\begin{equation}\label{eq:metric}
 |\nabla_{x'}\phi|^2=|\partial_r\phi|^2+r^{-2}|\partial_\varphi\phi|^2
 =\bigl(r\,\partial_\sigma\phi-(\lambda+Kr)\partial_z\phi\bigr)^2
    +\frac{(\partial_z\phi)^2}{(2\pi r)^2} .
\end{equation}

\subsection{The evolution of a material line}\label{sec:curve}

\begin{figure}[htbp]
\centering
\begin{tikzpicture}[x=1.16cm,y=1.16cm,line cap=round,line join=round]
  \pgfmathsetmacro{\outerdraw}{sqrt(13)/2}
  \newcommand{\drawspiralpart}[4]{%
    \draw[#1] plot[domain=#2:#3,samples=#4,variable=\w]
      ({\outerdraw*(.32+.68*\w/720)*cos(\w)},
       {\outerdraw*(.32+.68*\w/720)*sin(\w)});}
  \node[font=\small] at (0,2.20) {After the initial winding};
  \node[font=\small] at (5.65,2.20) {After the first pair of shears};
  \drawspiralpart{black!70,line width=.70pt}{0}{720}{241}
  \drawspiralpart{black,line width=3pt}{345.0000}{389.9610}{32}
  \begin{scope}[shift={(5.65,0)}]
    \drawspiralpart{black!55,dash pattern=on 2.6pt off 2.1pt,line width=.56pt}{0}{720}{241}
    \drawspiralpart{black,line width=3pt}{322.4312}{412.3542}{60}
  \end{scope}
\end{tikzpicture}
\caption[The tracked subarc before and after the first pair of shears.]{For illustration, we take $K=16$. The bold arcs track the same particles near $r=\sqrt3$ in the horizontal plane. The dashed curve is the original spiral.}
\label{fig:material}
\end{figure}
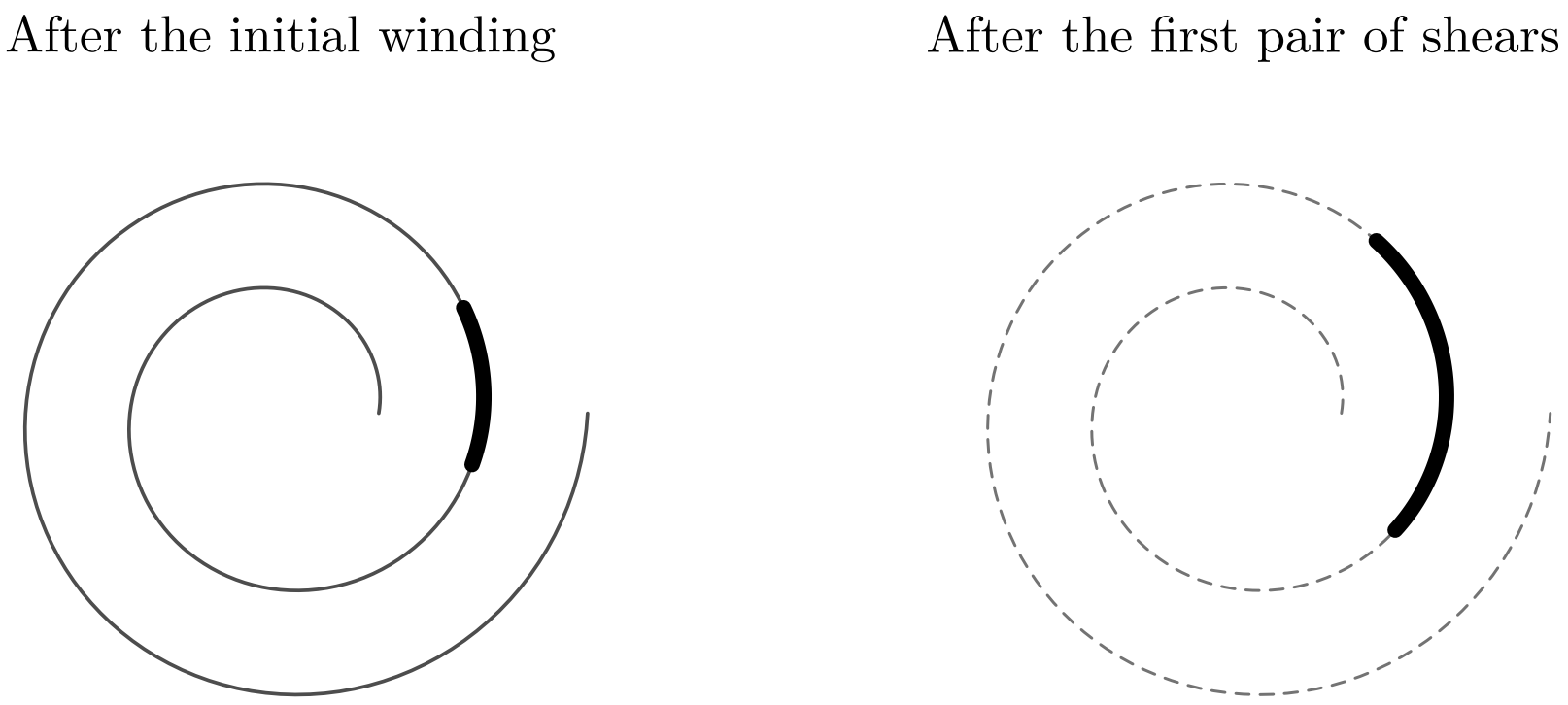

\begin{proof}[Proof of Proposition~\ref{prop:geometry}]

Fix $m\ge0$. On $[0,t_0]$, we use the winding velocity of Section~\ref{sec:preparation}. For $i=0,\ldots,m$, apply $R$ and $T$ on the first and second halves of $[t_i,t_{i+1}]$, respectively. Set $v^m=0$ on $[t_{m+1},T_*]$. Each time factor $\eta_I$ vanishes near the endpoints of $I$, so $v^m$ is smooth and $v^m(0)=0$. By Sections~\ref{sec:preparation} and~\ref{sec:shears}, $\operatorname{div}_{x'}v^m=0$ and $\operatorname{supp}v^m(t)\subset\Omega$.

During the initial winding, each marked particle keeps its radius $r$. By time $t_0$ it has made $\lambda(r-a)$ turns. The number of turns increases from zero to two as $r$ varies from $a$ to $b$. Thus the marked segment becomes the spiral $\Psi=0$. By \eqref{eq:radial-margin}, the marked particles stay where $\zeta=1$, so \eqref{eq:shears} applies in every pair.

At the start of pair $i$, take a marked particle on this spiral near $r=\sqrt3$, where $\sigma=3/2$ (Figure~\ref{fig:material}). Write its radial coordinate as $\sigma=3/2+d$. Set $q=2^i$, $p=Kq$, $h=h_i$, and take $|pd|<1/8$. Then we obtain $K(3/2+d)+z=0$ in $\mathbb R/\mathbb Z$.

The first shear $R$ keeps $z$ fixed, but the radius can change. Since $K$ is even, $3p/2$ is an integer. By periodicity and $S_h(s)=s$ for $|s|\le1/4$, we obtain $S_h(qz)=S_h(-pd)=-pd$. Substituting $S_h(qz)=-pd$ into \eqref{eq:shears} gives

\[
 R(3/2+d,z)=\left(3/2+d-\frac{-pd}{p},z\right)=(3/2+2d,z).
\]

Obviously, a particle moves outward if $d>0$ and moves inward if $d<0$. The angular difference from the spiral is now $\Psi=Kd$ turns.

The second shear $T$ keeps this radius fixed. Since $|2pd|<1/4$, we have $S_h(p(3/2+2d))=S_h(2pd)=2pd$. The angular change is therefore $-2pd/(2q)=-Kd$ turns. At the end of the pair,

\[
 \begin{aligned}
 T(3/2+2d,z)&=(3/2+2d,z-Kd),\\
 \Psi&=K(3/2+2d)+(z-Kd)=0 \qquad\text{in }\mathbb R/\mathbb Z .
 \end{aligned}
\]

Thus the particle returns to the original spiral at $\sigma=3/2+2d$. This coordinate belongs to $J$ since $|2d|<1/(4p)\le1/8$. The particle with $d=0$ stays still during both shears.

At $t_0$ the marked particles cover the original spiral for $|d|\le1/8$. For the first pair we need only $|d|<1/(8K)$. To follow the same particles through all $m+1$ pairs, take $|d|<1/(8K\,4^m)$ at $t_0$. For each $i=0,\ldots,m$, we have

\[
 (K2^i)\,|2^i d|=K\,4^i|d|\le K\,4^m|d|<\tfrac18 .
\]

For the marked particle with $\sigma(t_0)=3/2+d$, each pair of shears doubles the difference $\sigma-3/2$. The particle starts pair $i$ on the original spiral at $\sigma(t_i)=3/2+2^i d$ and returns to the same spiral at $\sigma(t_{i+1})=3/2+2^{i+1}d$ when the pair ends. As $d$ varies over the chosen interval, the radii after each pair fill an open interval in $(a,b)$. The marked particles thus cover a subarc of positive length.

\end{proof}

\FloatBarrier

\section{Scalar comparison on one shear interval}\label{sec:half}

Fix one shear of pair $n$, with $q=2^n$ and $p=Kq$. We measure time from the start of the shear and use the interval $[0,s]$, where $s$ is the duration of the shear. Denote its velocity by $v$. We compare the scalar at the start with $\phi_{\rm in}$ and at the end with $\phi_{\rm out}$. Using the functions of Section~\ref{sec:shears}, set

\[
 (\phi_{\rm in},\phi_{\rm out})=
 \begin{cases}
   (\phi_{p,q},\phi_{p,2q}), & R,\\
   (\phi_{p,2q},\phi_{2p,2q}), & T .
 \end{cases}
\]

For example, for the first radial shear we use $\phi_{K,1}$ and $\phi_{K,2}$. The angular shear then begins with the same comparison profile $\phi_{K,2}$ and ends with $\phi_{2K,2}$.

\begin{lemma}\label{lemma:half}

Fix a transition width $0<h<1/64$ and a diffusivity $\mu>0$. For any smooth initial datum, let $\theta$ be the solution on $[0,s]\times\T^2$ of

\[
 \partial_t\theta+v\cdot\nabla_{x'}\theta=\mu\Delta_{x'}\theta .
\]

For the error estimate, set

\[
 C_\chi=4\bigl(\|\chi\|_\infty+\|\chi^{\prime}\|_\infty\bigr).
\]

There is a constant $C_K$ for the diffusion error such that

\begin{equation}\label{eq:halfstep}
 \|\theta(s)-\phi_{\rm out}\|_2
 \le \|\theta(0)-\phi_{\rm in}\|_2+C_\chi\sqrt h+C_K(\mu s q^2)^{1/4}
       +\begin{cases}C_\chi/p & (R),\\ 0 & (T)\end{cases}.
\end{equation}

$C_K$ depends only on $K$, the annulus, $\lambda$, the $C^1$ norms of $\chi$ and $A$, and $\|\eta\|_\infty$. Both constants are independent of $n,s,h$ and $\mu$.

\end{lemma}

\begin{proof}

\textbf{Step 1. The transported scalar.} To estimate the error in each period, we use $M_\gamma$ from Section~\ref{sec:shears} with $2h\le\gamma<1/8$. Its definition on $[-1/2,1/2]$ gives

\[
 \begin{aligned}
 M_\gamma(\xi)&=1 && (|\xi|\le\tfrac12-2\gamma),\\
 M_\gamma(\xi)&=0 && (\tfrac12-\gamma\le|\xi|\le\tfrac12),\\
 |M_\gamma'(\xi)|&=\frac1\gamma
 \eta\!\left(\frac{1/2-|\xi|}{\gamma}-1\right)
 \le\frac{\|\eta\|_\infty}{\gamma} && (0<|\xi|<\tfrac12).
 \end{aligned}
\]

Clearly, $M_\gamma^{\prime}(0)=0$. Thus $M_\gamma=1$ except possibly within $2\gamma$ of either end of the period. For the two shears, set

\[
 M(\sigma,z)=
 \begin{cases}
 M_\gamma(qz) & (R),\\
 M_\gamma(p\sigma) & (T).
 \end{cases}
\]

For $\phi\in\{\phi_{\rm in},\phi_{\rm out}\}$, define $M\phi$ on $\T^2$ by

\[
 (M\phi)(x')=
 \begin{cases}
 M(\sigma,z)\,\phi(x'), & x'\in\Omega,\\
 0, & x'\in\T^2\setminus\Omega .
 \end{cases}
\]

Since $dx^{\prime}=d\sigma\,dz/(2\pi)$ and $M_\gamma$ has period one,

\[
 \begin{aligned}
 2\pi\,|\{M\ne1\}\cap\Omega|
 &=\begin{cases}
 \displaystyle\frac1q\big|\{\xi\in[0,q]: M_\gamma(\xi)\ne1\}\big| &(R),\\[4pt]
 \displaystyle\frac1p\big|\{\xi\in[p,2p]: M_\gamma(\xi)\ne1\}\big| &(T)
 \end{cases}\\
 &=\big|\{\xi\in[-1/2,1/2]:M_\gamma(\xi)\ne1\}\big|
 \le 2\cdot 2\gamma .
 \end{aligned}
\]

Here the substitutions are $\xi=qz$ in $R$ and $\xi=p\sigma$ in $T$. Since $0\le M\le1$ and $|\phi|\le\|\chi\|_\infty$,

\[
 \|(1-M)\phi\|_2^2
 =\int_{\Omega\cap\{M\ne1\}} |1-M|^2|\phi|^2\,dx'
 \le\|\chi\|_\infty^2\,|\{M\ne1\}\cap\Omega| .
\]

Taking square roots gives

\begin{equation}\label{eq:halfcut}
 |\{M\ne1\}\cap\Omega|\le\frac{4\gamma}{2\pi},\qquad
 \|(1-M)\phi\|_2\le\frac{2\|\chi\|_\infty}{\sqrt{2\pi}}\sqrt\gamma .
\end{equation}

Let $\rho$ solve the transport equation with initial value $M\phi_{\rm in}$,

\[
 \partial_t\rho+v\cdot\nabla_{x'}\rho=0,\qquad \rho(0)=M\phi_{\rm in}.
\]

For $0\le t\le s$, set $\tau=\tau_I(t_I+t)$. For either shear, \eqref{eq:ham} and $|S_h|\le1/2$ give

\[
 |\sigma(t)-\sigma(0)|
 \le \frac1{2p}\int_{t_I}^{t_I+t}\eta_I(u)\,du
 =\frac{\tau}{2p}\le\frac1{32}<\frac1{16} .
\]

Take a particle starting in $\operatorname{supp}(M\phi_{\rm in})$. Since $\operatorname{supp}\chi\subset(5/4,7/4)$ and $|\sigma(t)-\sigma(0)|<1/16$, we have

\[
 \frac{19}{16}=\frac54-\frac1{16}<\sigma(t)
 <\frac74+\frac1{16}=\frac{29}{16},\qquad 0\le t\le s.
\]

Thus $\zeta=1$ and $\zeta^{\prime}=0$ along this trajectory. By \eqref{eq:ham}, $z$ is constant during $R$, and $\sigma$ is constant during $T$. Write $(\sigma,z)$ for the position at time $t$. Integrating from $0$ to $t$ gives

\[
 (\sigma(0),z(0))=
 \begin{cases}
 (\sigma+\tau S_h(qz)/p,\ z) &(R),\\
 (\sigma,\ z+\tau S_h(p\sigma)/(2q)) &(T).
 \end{cases}
\]

The transport equation preserves $\rho$ along the trajectory. In both shears, $M$ depends only on the fixed coordinate. Hence

\[
 \rho(t,\sigma,z)=(M\phi_{\rm in})(\sigma(0),z(0))
 =M(\sigma,z)\,\phi_{\rm in}(\sigma(0),z(0)).
\]

Substituting the initial position into $\phi_{p,q}$ for $R$ and into $\phi_{p,2q}$ for $T$, we obtain

\[
 \rho=
 \begin{cases}
 M\,\chi\bigl(\sigma+\tau S_h(qz)/p\bigr)
      A\bigl(p\sigma+qz+\tau S_h(qz)\bigr) & (R),\\
 M\,\chi(\sigma)A\bigl(p\sigma+2qz+\tau S_h(p\sigma)\bigr) & (T).
 \end{cases}
\]

If the right-hand side is nonzero, the initial value of $\sigma$ lies in $\operatorname{supp}\chi$. Since $5/4<\sigma(0)<7/4$ and the radial displacement is at most $1/32<1/16$, the trajectory stays in $19/16<\sigma<29/16$, where $\zeta=1$ and $\zeta^{\prime}=0$. Integrating \eqref{eq:ham} forward from the initial coordinates gives

\[
 \begin{aligned}
 R:\quad \sigma(t)&=\left(\sigma+\frac{\tau}{p}S_h(qz)\right)
                       -\frac{\tau}{p}S_h(qz)=\sigma,\qquad z(t)=z,\\
 T:\quad z(t)&=\left(z+\frac{\tau}{2q}S_h(p\sigma)\right)
                       -\frac{\tau}{2q}S_h(p\sigma)=z,\qquad \sigma(t)=\sigma.
 \end{aligned}
\]

In each case the particle reaches $(\sigma,z)$ at time $t$. Thus the formula holds on $\Omega$.

Where $M\ne0$, use the following property of $M_\gamma$. On $[-1/2,1/2]$, $M_\gamma(\xi)\ne0$ implies $|\xi|<1/2-\gamma\le1/2-h$, where $S_h(\xi)=\xi$. By periodicity, if $M_\gamma(\xi)\ne0$, then $|\xi-k|<1/2-\gamma$ for some $k\in\mathbb Z$. For this $k$, $S_h(\xi)=S_h(\xi-k)=\xi-k$ and $S_h^{\prime}(\xi)=1$. Thus, wherever $M\ne0$, we have $S_h^{\prime}(qz)=1$ in $R$ and $S_h^{\prime}(p\sigma)=1$ in $T$.

In $R$, the partial derivatives of the argument of $\chi$ are

\[
 \begin{aligned}
 \partial_\sigma\!\left(\sigma+\frac{\tau}{p}S_h(qz)\right) &=1,\\
 \partial_z\!\left(\sigma+\frac{\tau}{p}S_h(qz)\right) &=\frac{\tau q}{p}S_h^{\prime}(qz)=\frac{\tau q}{p}.
 \end{aligned}
\]

For example, in $R$ the chain rule gives the following derivatives of the argument of $A$.

\[
 \begin{aligned}
 \partial_\sigma\bigl(p\sigma+qz+\tau S_h(qz)\bigr) &=p,\\
 \partial_z\bigl(p\sigma+qz+\tau S_h(qz)\bigr) &=q+\tau qS_h^{\prime}(qz)=(1+\tau)q.
 \end{aligned}
\]

In $T$, the argument of $\chi$ is $\sigma$, so

\[
 \begin{aligned}
 \partial_\sigma\,\sigma &=1,\\
 \partial_z\,\sigma &=0.
 \end{aligned}
\]

For the argument of $A$ in $T$, we obtain

\[
 \begin{aligned}
 \partial_\sigma\bigl(p\sigma+2qz+\tau S_h(p\sigma)\bigr) &=p+\tau pS_h^{\prime}(p\sigma)=(1+\tau)p,\\
 \partial_z\bigl(p\sigma+2qz+\tau S_h(p\sigma)\bigr) &=2q.
 \end{aligned}
\]

Differentiating the two formulas for $M$ gives

\[
 (\partial_\sigma M,\partial_z M)=
 \begin{cases}
  (0,qM_\gamma'(qz)) & (R),\\
  (pM_\gamma'(p\sigma),0) & (T).
 \end{cases}
\]

In each case one derivative is identically zero. Both vanish where $M=1$. Using $q\le p$, \eqref{eq:halfcut} and the bound for the derivative of $M_\gamma$, we obtain

\[
 \begin{aligned}
 &\|\partial_\sigma M\|_{L^2(\Omega)}+\|\partial_z M\|_{L^2(\Omega)}\\
 &\quad=\left(\int_{\Omega\cap\{M\ne1\}}
           (|\partial_\sigma M|^2+|\partial_z M|^2)\,dx'\right)^{1/2}\\
 &\quad\le p\|M_\gamma'\|_\infty\,|\Omega\cap\{M\ne1\}|^{1/2}
 \le \frac{p\|\eta\|_\infty}{\gamma}\left(\frac{4\gamma}{2\pi}\right)^{1/2}.
 \end{aligned}
\]

Using the formulas above, we obtain

\[
 \begin{aligned}
 \|\partial_\sigma\rho\|_2+\|\partial_z\rho\|_2
 &\le \|\chi\|_\infty\|A\|_\infty
       (\|\partial_\sigma M\|_2+\|\partial_z M\|_2)\\
 &\quad+|\Omega|^{1/2}(1+q/p)\|\chi'\|_\infty\|A\|_\infty\\
 &\quad+2|\Omega|^{1/2}(p+q)\|\chi\|_\infty\|A'\|_\infty\\
 &\le C\left(\frac p{\sqrt\gamma}+p+q\right).
 \end{aligned}
\]

The product rule, \eqref{eq:metric} and $p=Kq$ now yield

\begin{equation}\label{eq:halfgrad}
 \sup_{0\le t\le s}\|\nabla_{x'}\rho(t)\|_2
 \le C(1+K)\left(p+q+\frac{p}{\sqrt\gamma}\right)
 \le \frac{C_Kq}{\sqrt\gamma}.
\end{equation}

\textbf{Step 2. Comparison at the end of the shear.} Apply \eqref{eq:prep} on $[0,s]$ with $\rho(0)=M\phi_{\rm in}$. Taking square roots and using \eqref{eq:halfgrad}, we obtain

\[
 \begin{aligned}
 \|\theta(s)-\rho(s)\|_2
 &\le\|\theta(0)-M\phi_{\rm in}\|_2
          +\left(\mu\int_0^s\|\nabla_{x'}\rho(t)\|_2^2\,dt\right)^{1/2}\\
 &\le\|\theta(0)-\phi_{\rm in}\|_2+\|(1-M)\phi_{\rm in}\|_2
           +C_Kq\sqrt{\mu s/\gamma}.
 \end{aligned}
\]

At $t=s$, $\tau=1$. Where $M\ne0$, the input to $A$ in $\rho(s)$ differs from that in $\phi_{\rm out}$ by the integer $S_h(qz)-qz$ for $R$ and $S_h(p\sigma)-p\sigma$ for $T$. Since $A$ has period one, on $\Omega$ we obtain

\[
 \rho(s)-M\phi_{\rm out}=
 \begin{cases}
 M\bigl[\chi(\sigma+S_h(qz)/p)-\chi(\sigma)\bigr]A(p\sigma+2qz) &(R),\\[2pt]
 0 &(T).
 \end{cases}
\]

For $R$, use the mean value theorem and $|S_h|\le1/2$ to get

\[
 \left|\chi\!\left(\sigma+\frac{S_h(qz)}p\right)-\chi(\sigma)\right|
 \le\|\chi'\|_\infty\frac{|S_h(qz)|}{p}
 \le\frac{\|\chi'\|_\infty}{2p}.
\]

Since $|M A|\le1$, taking the $L^2$ norm gives

\[
 \|\rho(s)-M\phi_{\rm out}\|_2
 \le\begin{cases}
 |\Omega|^{1/2}\|\chi'\|_\infty/(2p) & (R),\\
 0 & (T).
 \end{cases}
\]

Substituting the estimate for $\theta(s)-\rho(s)$ into the triangle inequality gives

\[
 \begin{aligned}
 \|\theta(s)-\phi_{\rm out}\|_2
 &\le \|\theta(s)-\rho(s)\|_2+\|\rho(s)-M\phi_{\rm out}\|_2
             +\|(1-M)\phi_{\rm out}\|_2\\
 &\le \|\theta(0)-\phi_{\rm in}\|_2+C_Kq\sqrt{\mu s/\gamma}
             +\|\rho(s)-M\phi_{\rm out}\|_2\\
 &\qquad +\|(1-M)\phi_{\rm in}\|_2+\|(1-M)\phi_{\rm out}\|_2 .
 \end{aligned}
\]

Applying \eqref{eq:halfcut} to both profiles gives

\[
 \begin{aligned}
 \|(1-M)\phi_{\rm in}\|_2+\|(1-M)\phi_{\rm out}\|_2
 &\le 2\,\frac{2\|\chi\|_\infty}{\sqrt{2\pi}}\sqrt\gamma
 \le \frac{C_\chi}{\sqrt 2}\sqrt\gamma .
 \end{aligned}
\]

Since $|\Omega|^{1/2}\|\chi^{\prime}\|_\infty/2=\|\chi^{\prime}\|_\infty/(2\sqrt{2\pi})\le C_\chi$, substitution gives

\[
 \begin{aligned}
 \|\theta(s)-\phi_{\rm out}\|_2
 &\le\|\theta(0)-\phi_{\rm in}\|_2
    +\frac{C_\chi}{\sqrt 2}\sqrt\gamma+C_Kq\sqrt{\mu s/\gamma}
       +\begin{cases}C_\chi/p & (R),\\0 & (T).\end{cases}
 \end{aligned}
\]

If $q\sqrt{\mu s}<1/8$, take

\[
 \gamma=\max\{2h,q\sqrt{\mu s}\}.
\]

We have $2h\le\gamma<1/8$ since $h<1/64$. Also, $\sqrt{q\sqrt{\mu s}}=(\mu s q^2)^{1/4}$. Hence

\[
 \begin{aligned}
 \sqrt\gamma
 &=\max\{\sqrt{2h},(\mu s q^2)^{1/4}\}
 \le\sqrt{2h}+(\mu s q^2)^{1/4},\\
 q\sqrt{\mu s/\gamma}
 &=\frac{q\sqrt{\mu s}}{\sqrt{\gamma}}
 \le\frac{q\sqrt{\mu s}}{\sqrt{q\sqrt{\mu s}}}=(\mu s q^2)^{1/4}.
 \end{aligned}
\]

Substituting $p=Kq$ in \eqref{eq:halfgrad} and using $\sqrt\gamma<1$, we obtain

\[
 \begin{aligned}
 \sup_{0\le t\le s}\|\nabla_{x'}\rho(t)\|_2
 &\le C(1+K)\frac{q}{\sqrt\gamma}
             \bigl((K+1)\sqrt\gamma+K\bigr)\\
 &\le C(1+K)(2K+1)\frac{q}{\sqrt\gamma}.
 \end{aligned}
\]

Using this gradient bound in the diffusion term and substituting the two estimates for $\gamma$ gives

\[
 \begin{aligned}
 &\|\theta(s)-\phi_{\rm out}\|_2-\|\theta(0)-\phi_{\rm in}\|_2\\
 &\quad\le \frac{C_\chi}{\sqrt2}\bigl(\sqrt{2h}+(\mu s q^2)^{1/4}\bigr)
        + C(1+K)(2K+1)(\mu s q^2)^{1/4}
        +\begin{cases}C_\chi/p &(R),\\0 &(T)\end{cases}\\
 &\quad= C_\chi\sqrt h
        +\left(\frac{C_\chi}{\sqrt2}+C(1+K)(2K+1)\right)(\mu s q^2)^{1/4}
        +\begin{cases}C_\chi/p &(R),\\0 &(T).\end{cases}
 \end{aligned}
\]

If $q\sqrt{\mu s}\ge1/8$, the energy identity and the unit norms of both profiles give

\[
 \begin{aligned}
 \|\theta(s)\|_2^2
 &=\|\theta(0)\|_2^2-2\mu\int_0^s\|\nabla_{x'}\theta(t)\|_2^2\,dt
 \le\|\theta(0)\|_2^2,\\
 \|\theta(s)-\phi_{\rm out}\|_2
 &\le \|\theta(s)\|_2+\|\phi_{\rm out}\|_2\\
 &\le \|\theta(0)-\phi_{\rm in}\|_2+\|\phi_{\rm in}\|_2+\|\phi_{\rm out}\|_2\\
 &=\|\theta(0)-\phi_{\rm in}\|_2+2\\
 &\le \|\theta(0)-\phi_{\rm in}\|_2+2\sqrt8\,(\mu s q^2)^{1/4}.
 \end{aligned}
\]

In the first case, the two terms proportional to $(\mu s q^2)^{1/4}$ have coefficients $C_\chi/\sqrt2$ and $C(1+K)(2K+1)$. The second estimate has coefficient $4\sqrt2$. Taking the larger coefficient then yields

\[
 C_K=\max\left\{\frac{C_\chi}{\sqrt2}+C(1+K)(2K+1),\ 4\sqrt2\right\}.
\]

By choosing $C_K$ as above, we obtain \eqref{eq:halfstep} in the first case. In the second case, the estimate above already gives

\[
 \|\theta(s)-\phi_{\rm out}\|_2
 \le\|\theta(0)-\phi_{\rm in}\|_2+C_K(\mu s q^2)^{1/4}.
\]

For $R$, the right-hand side of \eqref{eq:halfstep} contains the additional terms $C_\chi\sqrt h+C_\chi/p$. For $T$, it contains only the additional term $C_\chi\sqrt h$. These terms are nonnegative, so \eqref{eq:halfstep} follows for either shear. This proves the lemma.

\end{proof}

\section{Anomalous dissipation for the scalar}\label{sec:scalar-dissipation}

We use Lemma~\ref{lemma:half} to compare the scalar with the spiral profile until the horizontal flow stops. We then estimate the diffusion of this profile on the remaining time interval.

\begin{proposition}[Scalar dissipation in the subannulus]\label{prop:scalar}
Let $\theta_0$ be the smooth initial datum in \eqref{eq:initial-profile}, with $\int_{\T^2}\theta_0\,dx'=0$, $\|\theta_0\|_2=1$. The velocities in Proposition~\ref{prop:geometry} can be chosen together with a sequence $\mu_m\downarrow0$ so that the solutions $\theta^m$ of \eqref{eq:scalar} with $\theta^m(0)=\theta_0$ satisfy

\[
 \liminf_{m\to\infty}\mu_m\int_0^{T_*}\!\int_{\Omega'}|\nabla\theta^m|^2\,dx'\,dt\ge\tfrac38 .
\]
\end{proposition}

\subsection{Time intervals and viscosity}

Set $t_0\in(0,T_*)$ to be the time when the initial winding ends. For each $n\ge0$, set the duration of each shear in pair $n$ to $s_n$, where

\begin{equation}\label{eq:timetable}
 t_n=T_*-\frac{T_*-t_0}{n+1},\qquad 2s_n=t_{n+1}-t_n=\frac{T_*-t_0}{(n+1)(n+2)}.
\end{equation}

Thus $t_n\uparrow T_*$ and the shear intervals fill $[t_0,T_*)$. Following \cite[Sec.~5.1]{BDL}, stop $v^m$ after pair $m$. The time for diffusion is

\begin{equation}\label{eq:window}
 R_m=T_*-t_{m+1}=\frac{T_*-t_0}{m+2}.
\end{equation}

After each pair, the angular frequency of the comparison profile doubles, so $q_n=2^n$. By \eqref{eq:timetable} and \eqref{eq:window}, the duration $s_m$ of the last shear and the time $R_m$ left for diffusion satisfy

\[
 s_m=\frac{T_*-t_0}{2(m+1)(m+2)},\qquad \frac{R_m}{s_m}=2(m+1).
\]

As in \cite[Sec.~2.6]{CCS}, we choose the viscosity to allow little diffusion during the last shear and substantial decay of the profile after the flow stops.

Set $\mu_m q_m^2=(m+1)^{3/2}$. Since $3/2<2$, we have $\mu_m s_m q_m^2\to0$, which controls the diffusion error during the last shear. Since $3/2>1$ and $q_{m+1}=2q_m$, we have $\mu_m R_m q_{m+1}^2\to\infty$, which gives decay of the final spiral profile after the flow stops. We take

\begin{equation}\label{eq:viscosities}
 \mu_m=(m+1)^{3/2}4^{-m},\qquad
 h_n=\left(\frac{2s_n}{K(T_*-t_0)}\right)^2=\frac{1}{K^2(n+1)^2(n+2)^2}.
\end{equation}

The error term from the transition intervals in Lemma~\ref{lemma:half} is bounded by $C_\chi\sqrt{h_n}$. By \eqref{eq:viscosities},

\[
 \sum_{n=0}^m\sqrt{h_n}=\frac{t_{m+1}-t_0}{K(T_*-t_0)}=\frac{1}{K}\!\left(1-\frac1{m+2}\right)\le\frac1K.
\]

Since $K\ge16$, every width satisfies $0<h_n\le 1/(4K^2)<1/64$, as we need.

By \eqref{eq:timetable}--\eqref{eq:viscosities} and $m+2\le2(m+1)$,

\[
 \begin{aligned}
 \mu_m s_m q_m^2&=\frac{(T_*-t_0)\sqrt{m+1}}{2(m+2)}\longrightarrow0,\\
 \mu_m R_m q_{m+1}^2&=\frac{4(T_*-t_0)(m+1)^{3/2}}{m+2}
 \ge 2(T_*-t_0)\sqrt{m+1}\longrightarrow\infty.
 \end{aligned}
\]

Also $\mu_{m+1}/\mu_m=\tfrac14(1+1/(m+1))^{3/2}\le2^{-1/2}$, so $\mu_m\downarrow0$.

\subsection{The scalar error when the flow stops}

For the velocities $v^m$ constructed in Section~\ref{sec:geometry} and the same initial datum $\theta_0=\chi(\sigma)A(y)$ in \eqref{eq:initial-profile}, equation~\eqref{eq:scalar} has a unique smooth solution by the periodic version of the parabolic estimates in \cite[Secs.~7.1.2--7.1.3]{EV}. The scalar energy identity with $\operatorname{div}v^m=0$ and $\|\theta_0\|_2=1$ yields

\[
 \|\theta^m(t)\|_2^2+2\mu_m\int_0^t\|\nabla\theta^m(s)\|_2^2\,ds=1,
 \qquad 0\le t\le T_* .
\]

At $t_{m+1}$, set $\phi^m=\phi_{Kq_{m+1},q_{m+1}}$ on $\T^2$. By definition, $\phi^m=\chi(\sigma)A(q_{m+1}\Psi)$ on $\Omega$ and $\phi^m=0$ outside $\Omega$. On each shear interval take a second viscous solution with initial value $\phi_{\rm in}$.

On $I=(t_I,t_I+s_n)$ denote this solution by $\widetilde\theta$ and set $w=\theta^m-\widetilde\theta$. Both solutions have velocity $v^m$ and diffusivity $\mu_m$, so

\[
 \partial_t w+v^m\cdot\nabla_{x'}w=\mu_m\Delta_{x'}w,
 \qquad w(t_I)=\theta^m(t_I)-\phi_{\rm in}.
\]

Taking the inner product with $2w$, using $\operatorname{div}_{x'}v^m=0$ and integrating in time gives

\[
 \|w(t)\|_2^2+2\mu_m\int_{t_I}^{t}\|\nabla_{x'}w(\tau)\|_2^2\,d\tau
 =\|\theta^m(t_I)-\phi_{\rm in}\|_2^2,\qquad t_I\le t\le t_I+s_n .
\]

The integral is nonnegative. At the end of the shear the triangle inequality therefore gives

\[
 \begin{aligned}
 \|\theta^m(t_I+s_n)-\phi_{\rm out}\|_2
 &\le \|w(t_I+s_n)\|_2+\|\widetilde\theta(t_I+s_n)-\phi_{\rm out}\|_2\\
 &\le \|\theta^m(t_I)-\phi_{\rm in}\|_2
       +\|\widetilde\theta(t_I+s_n)-\phi_{\rm out}\|_2 .
 \end{aligned}
\]

Applying \eqref{eq:halfstep} to both shears in each pair, summing over $n=0,\ldots,m$, and using \eqref{eq:prep} for the initial winding, we obtain

\begin{align}
 \|\theta^m(t_{m+1})-\phi^m\|_2
 &\le C\sqrt{\mu_m}
 +2C_\chi\sum_{n=0}^m\sqrt{h_n}
 +\frac{C_\chi}{K}\sum_{n=0}^m2^{-n}\notag\\
 &\quad+2C_K\sum_{n=0}^m(\mu_m s_nq_n^2)^{1/4}.
 \label{eq:sumhalf}
\end{align}

In pair $n$, set $q=2^n$ and $p=Kq$ as in Section~\ref{sec:shears}. On the support of the transported profile, the first shear $R$ keeps $z$ fixed. A particle that ends this shear at $(\sigma,z)$ started at $(\sigma+S_{h_n}(qz)/p,z)$. Since $|S_{h_n}|\le1/2$, the mean value theorem gives

\[
 \left|\chi\!\left(\sigma+\frac{S_{h_n}(qz)}{p}\right)-\chi(\sigma)\right|
 \le\|\chi'\|_\infty\frac{|S_{h_n}(qz)|}{p}
 \le\frac{\|\chi'\|_\infty}{2p}.
\]

The second shear $T$ keeps $\sigma$ fixed, so its radial factor is still $\chi(\sigma)$. Thus the displacement term $C_\chi/p=C_\chi/(K2^n)$ in \eqref{eq:halfstep} occurs once per pair.

For the viscous sum, \eqref{eq:timetable} and \eqref{eq:viscosities} give

\[
 \begin{aligned}
 (\mu_m s_n q_n^2)^{1/4}
 &=\left(\frac{T_*-t_0}{2}\right)^{1/4}
   \frac{(m+1)^{3/8}\,2^{-(m-n)/2}}{[(n+1)(n+2)]^{1/4}}\\
 &\le\left(\frac{T_*-t_0}{2}\right)^{1/4}
   (m+1)^{3/8}(n+1)^{-1/2}2^{-(m-n)/2}.
 \end{aligned}
\]

Set $j=m-n$. Since $(m+1)/(n+1)\le j+1$ and $\sqrt{j+1}\le j+1$,

\[
 \sum_{n=0}^m(n+1)^{-1/2}2^{-(m-n)/2}
 \le (m+1)^{-1/2}\sum_{j\ge0}(j+1)2^{-j/2}
 =\frac{(m+1)^{-1/2}}{(1-2^{-1/2})^2}.
\]

The viscous sum is bounded by $C_K(m+1)^{3/8-1/2}$. The other two sums in \eqref{eq:sumhalf} satisfy

\[
 \begin{aligned}
 2C_\chi\sum_{n=0}^m\sqrt{h_n}
 &=\frac{2C_\chi}{K}\!\left(1-\frac1{m+2}\right)\le\frac{2C_\chi}{K},\\
 \frac{C_\chi}{K}\sum_{n=0}^m 2^{-n}
 &=\frac{C_\chi}{K}\,\frac{1-2^{-(m+1)}}{1-1/2}
 \le\frac{2C_\chi}{K}.
 \end{aligned}
\]

Let $E_*$ denote the sum of these two upper bounds. Substituting these bounds into \eqref{eq:sumhalf} gives

\begin{equation}\label{eq:stoperror}
 \|\theta^m(t_{m+1})-\phi^m\|_2
 \le E_*+C\sqrt{\mu_m}+C_K(m+1)^{-1/8},\qquad
 E_*=\frac{4C_\chi}{K}.
\end{equation}

By Lemma~\ref{lemma:half}, $C_\chi$ depends only on the fixed function $\chi$. Choose one even integer $K\ge16$ with $K>16C_\chi$. The first inequality keeps the transported support in the region where $\zeta=1$ by \eqref{eq:radial-margin}. The second gives $E_*=4C_\chi/K<1/4$ in \eqref{eq:stoperror}.

In \eqref{eq:stoperror}, we absorb the factor $2((T_*-t_0)/2)^{1/4}/(1-2^{-1/2})^2$ into the constant from Lemma~\ref{lemma:half}. The resulting constant $C_K$ depends only on $K$, $T_*-t_0$, the annulus, the slope $\lambda$, the $C^1$ norms of $\chi$ and $A$, and $\|\eta\|_\infty$. It is independent of $m$, $n$, $h_n$ and $\mu_m$.

With $K$ fixed, $\mu_m\to0$ and $(m+1)^{-1/8}\to0$. Hence \eqref{eq:stoperror} gives

\[
 \limsup_{m\to\infty}\|\theta^m(t_{m+1})-\phi^m\|_2
 \le E_*<\frac14 .
\]

For all sufficiently large $m$, the $L^2$ error is below $1/4$ at $t_{m+1}$ and remains so under diffusion. By the triangle inequality, $\|\theta^m(T_*)\|_2\le1/2$ follows once the solution starting from $\phi^m$ has norm at most $1/4$ at $T_*$.

\subsection{Scalar diffusion after the last shear}

On $[t_{m+1},T_*]$, substituting $v^m=0$ into \eqref{eq:scalar} yields $\partial_t\theta^m=\mu_m\Delta\theta^m$. Cheskidov \cite[Sec.~4]{CH} estimated the scalar on such an interval by showing that its low Fourier modes were already small when the flow stopped and that diffusion damped its high modes. To estimate the decay of our spiral profile, we use $\partial_\varphi\cos(2\pi q\Psi)=-q\sin(2\pi q\Psi)$ to write the profile as the divergence of a vector field.

Write $(x')^\perp=(-x_2,x_1)$. Its divergence is zero and $(x')^\perp\cdot\nabla=\partial_\varphi$ on the annulus. As $\partial_\varphi(2\pi q\Psi)=q$ and $\chi$ is radial, for each integer $q\ge1$ we have

\begin{equation}\label{eq:angular-primitive}
 H_q=-\frac1q\chi(\sigma)\cos(2\pi q\Psi)(x')^\perp,\qquad
 \operatorname{div}H_q=\chi(\sigma)\sin(2\pi q\Psi),\qquad
 \|H_q\|_2\le\frac{2}{q}.
\end{equation}

Set $H_q=0$ outside $\Omega$.

At fixed $r$, the phase $2\pi q\Psi$ runs through $q$ full periods as $y$ runs from $0$ to $1$. Hence

\[
 \int_0^1\cos^2(2\pi q\Psi)\,dy
 =\int_0^1\sin^2(2\pi q\Psi)\,dy=\frac12 .
\]

Since $|(x')^\perp|=|x'|\le2$ on the support of $\chi$, the normalized area element and $\int_1^2\chi^2\,d\sigma=4\pi$ give

\[
 \begin{aligned}
 \|H_q\|_2^2
 &=\frac{1}{2\pi q^2}\int_1^2 \chi(\sigma)^2 r^2
           \int_0^1\cos^2(2\pi q\Psi)\,dy\,d\sigma\\
 &\le\frac{4}{2\pi q^2}\,\frac12\int_1^2\chi(\sigma)^2\,d\sigma
 =\frac4{q^2}.
 \end{aligned}
\]

For a smooth scalar or vector field $F$ on $\T^2$, let $e^{\tau\Delta}F$ be the solution at time $\tau>0$ of the diffusion equation with diffusivity one and initial value $F$. In Fourier variables, for $k\in\mathbb Z^2$,

\[
 \widehat F(k)=\int_{\T^2}F(x')e^{-ik\cdot x'}\,dx',\qquad
 \widehat{e^{\tau\Delta}F}(k)=e^{-\tau|k|^2}\widehat F(k).
\]

Integration by parts gives $\widehat{\operatorname{div}H_q}(k)=ik\cdot\widehat H_q(k)$. Parseval's identity and $|k\cdot\widehat H_q(k)|\le|k|\,|\widehat H_q(k)|$ give

\[
 \begin{aligned}
 \|e^{\tau\Delta}\operatorname{div}H_q\|_2^2
 &=\sum_{k\in\mathbb Z^2}e^{-2\tau|k|^2}|k\cdot\widehat H_q(k)|^2\\
 &\le\sum_{k\in\mathbb Z^2}|k|^2 e^{-2\tau|k|^2}|\widehat H_q(k)|^2
 \le \Bigl(\sup_{x\ge0}x e^{-2\tau x}\Bigr)\|H_q\|_2^2 .
 \end{aligned}
\]

The derivative of $xe^{-2\tau x}$ is $(1-2\tau x)e^{-2\tau x}$, so the maximum is attained at $x=1/(2\tau)$. Taking square roots and using \eqref{eq:angular-primitive}, we obtain

\[
 \sup_{x\ge0}xe^{-2\tau x}=\frac{1}{2e\tau},\qquad
 \|e^{\tau\Delta}\operatorname{div}H_q\|_2
 \le\frac{\|H_q\|_2}{\sqrt{2e\tau}}
 \le\frac{\sqrt{2/e}}{q\sqrt\tau}.
\]

Take $q=q_{m+1}$ and $\tau=\mu_mR_m$. Then $\phi^m=\operatorname{div}H_{q_{m+1}}$. By \eqref{eq:window} and \eqref{eq:viscosities}, $\mu_mR_m q_{m+1}^2\ge 2(T_*-t_0)\sqrt{m+1}$. Taking the inverse square root and using the bound above, we obtain

\begin{equation}\label{eq:diffusion-bound}
 \begin{aligned}
 \|e^{\mu_mR_m\Delta}\phi^m\|_2
 &=\|e^{\mu_mR_m\Delta}\operatorname{div}H_{q_{m+1}}\|_2
 \le\frac{\sqrt{2/e}}{q_{m+1}\sqrt{\mu_mR_m}}\\
 &\le \frac{\sqrt{2/e}}{\sqrt{2(T_*-t_0)}}(m+1)^{-1/4}
 =\frac{1}{\sqrt{e(T_*-t_0)}}(m+1)^{-1/4}.
 \end{aligned}
\end{equation}

\subsection{Dissipation in the subannulus}

We now complete the proof of Proposition~\ref{prop:scalar} by estimating the dissipation in $\Omega'$.

\begin{proof}[Proof of Proposition~\ref{prop:scalar}]

Since $\theta^m(T_*)=e^{\mu_mR_m\Delta}\theta^m(t_{m+1})$, we separate the error at $t_{m+1}$ from the spiral profile. By Parseval's identity, $e^{-\tau|k|^2}\le1$ gives $\|e^{\tau\Delta}F\|_2\le\|F\|_2$. Using \eqref{eq:stoperror} and \eqref{eq:diffusion-bound}, we obtain

\[
 \begin{aligned}
 \|\theta^m(T_*)\|_2
 &\le\|e^{\mu_mR_m\Delta}(\theta^m(t_{m+1})-\phi^m)\|_2+\|e^{\mu_mR_m\Delta}\phi^m\|_2\\
 &\le \|\theta^m(t_{m+1})-\phi^m\|_2+\|e^{\mu_mR_m\Delta}\phi^m\|_2\\
 &\le E_*+C\sqrt{\mu_m}+C_K(m+1)^{-1/8}+C(m+1)^{-1/4}.
 \end{aligned}
\]

Thus $\|\theta^m(T_*)\|_2\le1/2$ for all sufficiently large $m$, since $E_*<1/4$.

Choose $0\le P\le1$ in $C_c^\infty(\Omega')$, equal to one near the supports of $\zeta$ and $\chi$, and extend it by zero. Then $P\theta_0=\theta_0$ and $v^m\cdot\nabla P=0$.

Taking the $L^2(\T^2)$ inner product of \eqref{eq:scalar} with $2P\theta^m$ gives

\[
 \frac{d}{dt}\int_{\T^2}P|\theta^m|^2\,dx'
 +\int_{\T^2}P v^m\cdot\nabla_{x'}|\theta^m|^2\,dx'
 =2\mu_m\int_{\T^2}P\theta^m\Delta_{x'}\theta^m\,dx'.
\]

For the transport term, periodic integration by parts gives

\[
 \int_{\T^2}P v^m\cdot\nabla_{x'}|\theta^m|^2\,dx'
 =-\int_{\T^2}|\theta^m|^2
 (v^m\cdot\nabla_{x'}P+P\operatorname{div}_{x'}v^m)\,dx'=0 .
\]

For the diffusion term, integrating by parts twice gives

\[
 \begin{aligned}
 2\mu_m\int_{\T^2}P\theta^m\Delta_{x'}\theta^m\,dx'
 &=-2\mu_m\int_{\T^2}P|\nabla_{x'}\theta^m|^2\,dx'
   -2\mu_m\int_{\T^2}\theta^m\nabla_{x'}P\cdot\nabla_{x'}\theta^m\,dx'\\
 &=-2\mu_m\int_{\T^2}P|\nabla_{x'}\theta^m|^2\,dx'
   -\mu_m\int_{\T^2}\nabla_{x'}P\cdot\nabla_{x'}|\theta^m|^2\,dx'\\
 &=-2\mu_m\int_{\T^2}P|\nabla_{x'}\theta^m|^2\,dx'
   +\mu_m\int_{\T^2}(\Delta_{x'}P)|\theta^m|^2\,dx'.
 \end{aligned}
\]

Combining these identities yields

\[
 \frac{d}{dt}\int_{\T^2} P|\theta^m|^2\,dx'
 +2\mu_m\int_{\T^2}P|\nabla_{x'}\theta^m|^2\,dx'
 =\mu_m\int_{\T^2}(\Delta_{x'}P)|\theta^m|^2\,dx'.
\]

At time zero the integral of $P|\theta_0|^2$ equals one. At $T_*$, the bound $0\le P\le1$ gives

\[
 \int_{\T^2}P(x')|\theta^m(T_*,x')|^2\,dx'
 \le\int_{\T^2}|\theta^m(T_*,x')|^2\,dx'
 =\|\theta^m(T_*)\|_2^2 .
\]

Since $\|\theta^m(t)\|_2\le1$, time integration gives

\[
 2\mu_m\int_0^{T_*}\!\int_{\Omega'}|\nabla\theta^m|^2\,dx'\,dt
 \ge 1-\|\theta^m(T_*)\|_2^2-\mu_m T_*\|\Delta P\|_\infty .
\]

The cutoff error tends to zero since $P$ is fixed. The initial scalar energy is $1/2$, and the energy at $T_*$ is at most $1/8$ for all sufficiently large $m$. Taking the lower limit gives the lower bound

\begin{equation}\label{eq:scalardiss}
\liminf_{m\to\infty}\mu_m\int_0^{T_*}\!\int_{\Omega'}|\nabla\theta^m|^2\,dx'\,dt\ge\tfrac38 .
\end{equation}
\end{proof}

\section{Anomalous dissipation for the forced 3D Navier--Stokes equations}\label{sec:NS}

We estimate the horizontal forces and use the scalar dissipation bound to prove Theorem~\ref{thm:main}.

\subsection{Convergence of the forces}

We follow the argument in \cite[Sec.~5.1]{BDL} to prove the convergence of the forces. Take the velocities $v^m$ and viscosities $\mu_m$ constructed in Section~\ref{sec:scalar-dissipation}. Define
\begin{equation}\label{eq:nsforce}
 g^m=\partial_t v^m+(v^m\cdot\nabla)v^m-\mu_m\Delta v^m .
\end{equation}

\begin{proposition}\label{prop:force}
There exists a horizontal force $g$ such that, for every $0<\alpha<1$,
\[
 g^m\longrightarrow g\quad\text{in }C([0,T_*];C^\alpha(\T^2)).
\]
\end{proposition}

\begin{proof}[Proof]
Fix $0<\alpha<1$. On a shear interval $I$ of pair $n$, let $q=2^n$, $p=Kq$, $h=h_n$ and $v^m=\eta_I w$. By periodicity, it suffices to bound $Q_h$ on $[0,1]$. Its defining integral gives $\|Q_h\|_\infty\le\|S_h\|_\infty\le1/2$. Thus \eqref{eq:ham} and the velocity formula in Section~\ref{sec:shears} give $\|w\|_\infty\le C/q$. For $1\le j\le3$, the chain rule and $p=Kq$ give

\[
 \begin{aligned}
 \|\partial_z^j S_h(qz)\|_\infty
   &=q^j\|S_h^{(j)}\|_\infty\le C_j(q/h)^j,\\
 \|\partial_\sigma^j S_h(p\sigma)\|_\infty
   &\le p^j\|S_h^{(j)}\|_\infty\le C_j K^j(q/h)^j .
 \end{aligned}
\]

The same bounds hold for $Q_h$ since $Q_h^{(j)}=S_h^{(j-1)}$. Applying the product rule to the velocity formula in Section~\ref{sec:shears} gives, for $0\le\ell\le3$,

\begin{equation}\label{eq:ns-space}
 \|\nabla_{x'}^\ell w\|_\infty
 \le \frac{C}{q}\sum_{j=0}^{\ell}(q/h)^j
 \le C h^{-\ell}q^{\ell-1}.
\end{equation}

The constant is independent of $n,m$. By the choices in Section~\ref{sec:scalar-dissipation},

\[
 s_n^{-1}=\frac{2(n+1)(n+2)}{T_*-t_0}\le\frac{4(n+1)^2}{T_*-t_0},
 \qquad h_n^{-1}=K^2(n+1)^2(n+2)^2\le4K^2(n+1)^4 .
\]

Using $\|\eta_I^{(j)}\|_\infty\le Cs_n^{-1-j}$ for $j=0,1$ then gives

\begin{equation}\label{eq:ns-derivatives}
\begin{aligned}
 \|\nabla^\ell v^m\|_\infty
    &\le C(n+1)^{2+4\ell}q^{\ell-1} &&(0\le\ell\le3),\\
 \|\nabla^\ell\partial_t v^m\|_\infty
    &\le C(n+1)^{4+4\ell}q^{\ell-1} &&(0\le\ell\le1).
\end{aligned}
\end{equation}

Set $G^m=\partial_t v^m+(v^m\cdot\nabla)v^m$. By the product rule and \eqref{eq:ns-derivatives},

\[
\begin{aligned}
 \|G^m\|_\infty
 &\le \|\partial_t v^m\|_\infty+\|v^m\|_\infty\|\nabla v^m\|_\infty\\
 &\le C\left(\frac{(n+1)^4}{q}
       +\frac{(n+1)^2}{q}(n+1)^6\right)
 \le C\frac{(n+1)^8}{q},\\[2pt]
 \|\nabla G^m\|_\infty
 &\le \|\nabla\partial_t v^m\|_\infty
       +C\bigl(\|\nabla v^m\|_\infty^2+
                 \|v^m\|_\infty\|\nabla^2 v^m\|_\infty\bigr)\\
 &\le C\left((n+1)^8+(n+1)^{12}
          +\frac{(n+1)^2}{q}(n+1)^{10}q\right)
 \le C(n+1)^{12}.
\end{aligned}
\]

Taking $\ell=2$ in \eqref{eq:ns-derivatives} gives $\|\Delta v^m\|_\infty\le C(n+1)^{10}q$. For $0<\alpha<1$, we use

\[
 \|F\|_{C^\alpha}\le\|F\|_\infty+
        (2\|F\|_\infty)^{1-\alpha}\|\nabla F\|_\infty^\alpha .
\]

Using the bounds above and the case $\ell=3$ of \eqref{eq:ns-derivatives}, we obtain

\[
 \begin{aligned}
 \|G^m\|_{C^\alpha}
 &\le C\left[\frac{(n+1)^8}{q}
 +\left(\frac{(n+1)^8}{q}\right)^{1-\alpha}
       \bigl((n+1)^{12}\bigr)^\alpha\right]
 \le C(n+1)^{8+4\alpha}q^{\alpha-1},\\[3pt]
 \|\Delta v^m\|_{C^\alpha}
 &\le C\left[(n+1)^{10}q+
       \bigl((n+1)^{10}q\bigr)^{1-\alpha}
       \bigl((n+1)^{14}q^2\bigr)^\alpha\right]
 \le C(n+1)^{10+4\alpha}q^{1+\alpha}.
 \end{aligned}
\]

For $k>m$, the construction gives

\[
 v^m(t)=
 \begin{cases}
 v^k(t), & 0\le t\le t_{m+1},\\
 0, & t_{m+1}<t\le T_* .
 \end{cases}
\]

Only pairs $m+1,\ldots,k$ contribute to the difference. The bound for $G^k$ on each shear interval then gives

\[
 \sup_{t\in[0,T_*]}\|G^k(t)-G^m(t)\|_{C^\alpha}
 \le C\sup_{n>m}(n+1)^{8+4\alpha}2^{-(1-\alpha)n}
 \longrightarrow0 \qquad(m\to\infty).
\]

Thus $G^m$ converges to its uniform limit $g$ in $C([0,T_*];C^\alpha(\T^2))$. On $[0,t_0]$ the winding velocity is the same for every $m$, so $\mu_m\|\Delta v^m(t)\|_{C^\alpha}\le C\mu_m$. For the shears, $n\le m$, so \eqref{eq:viscosities} gives
\[
 \sup_{t\in[0,T_*]}\mu_m\|\Delta v^m(t)\|_{C^\alpha}
 \le C\mu_m+C(m+1)^{23/2+4\alpha}2^{-(1-\alpha)m}\longrightarrow0 .
\]
By \eqref{eq:nsforce}, the forces $g^m$ converge to the same $g$ for every $0<\alpha<1$.

\end{proof}

\subsection[Proof of the main theorem]{Proof of Theorem~\ref{thm:main}}

\begin{proof}[Proof]
Using the $(2+\tfrac12)$-dimensional construction of Bruè and De Lellis \cite[Sec.~3.1]{BDL} with the velocity $v^m$ and scalar $\theta^m$ from the proof of Proposition~\ref{prop:scalar}, we set
\[
 \begin{aligned}
 u^m(t,x',x_3)&=(v^m(t,x'),\theta^m(t,x')),\\
 f^m(t,x',x_3)&=(g^m(t,x'),0),\qquad p^m(t,x',x_3)=0 .
 \end{aligned}
\]
Since these fields are independent of $x_3$, we reduce \eqref{eq:NS} to
\[
 \begin{cases}
 \partial_t v^m+(v^m\cdot\nabla)v^m=\mu_m\Delta v^m+g^m,\\
 \operatorname{div}v^m=0,\\
 \partial_t\theta^m+v^m\cdot\nabla\theta^m=\mu_m\Delta\theta^m .
 \end{cases}
\]
We used the first equation to define $g^m$ in \eqref{eq:nsforce} and the third equation, \eqref{eq:scalar}, to prove Proposition~\ref{prop:scalar}. The divergence-free condition follows from Proposition~\ref{prop:geometry}. Since $v^m(0)=0$ and $\theta^m(0)=\theta_0$, the smooth solutions have the same initial datum $(0,0,\theta_0)$.

Let $g$ be the limit in Proposition~\ref{prop:force}. Define its extension by
\[
 f(t,x',x_3)=(g(t,x'),0)\qquad(x_3\in\T).
\]
Fix $t$ and $0<\alpha<1$. Let $L$ be any Hölder constant for $(g^m-g)(t)$. For $x=(x',x_3)$ and $y=(y',y_3)$ in $\T^3$, the forces depend only on the horizontal coordinates and $|x'-y'|\le |x-y|$. Hence
\[
 \begin{aligned}
 |(f^m-f)(t,x)-(f^m-f)(t,y)|
 &=|(g^m-g)(t,x')-(g^m-g)(t,y')|\\
 &\le L|x'-y'|^\alpha\le L|x-y|^\alpha .
 \end{aligned}
\]
This holds for every such $L$, and the two force differences have the same supremum norm. Taking the supremum in time and applying Proposition~\ref{prop:force} gives
\[
 \sup_{0\le t\le T_*}\|f^m(t)-f(t)\|_{C^\alpha(\T^3)}
 \le \sup_{0\le t\le T_*}\|g^m(t)-g(t)\|_{C^\alpha(\T^2)}
 \longrightarrow 0 \qquad (m\to\infty).
\]
This proves that $f^m\to f$ in $C([0,T_*];C^\alpha(\T^3))$.

Since $\partial_3u^m=0$ and $\int_{\T}1\,dx_3=1$,
\[
 \begin{aligned}
 \mu_m\int_0^{T_*}\!\int_{\Omega'\times\T}|\nabla u^m|^2\,dx'\,dx_3\,dt
 &=\mu_m\int_0^{T_*}\!\int_{\Omega'}
    \bigl(|\nabla v^m|^2+|\nabla\theta^m|^2\bigr)\,dx'\,dt\\
 &\ge \mu_m\int_0^{T_*}\!\int_{\Omega'}|\nabla\theta^m|^2\,dx'\,dt .
 \end{aligned}
\]
Taking the lower limit as $m\to\infty$ and applying Proposition~\ref{prop:scalar}, we obtain the lower bound $3/8$ in \eqref{eq:main-local}. Since $\Omega'\times\T\subset\T^3$, the three-dimensional solutions exhibit anomalous dissipation in the sense of \eqref{eq:ad-definition}. This completes the proof of Theorem~\ref{thm:main}.

\end{proof}

\noindent\textbf{Acknowledgments.} The author is deeply grateful to Professor Liqun Zhang for his guidance and support.

\vspace{1em}
\noindent\textbf{Declaration of Generative Artificial Intelligence (AI) Use.} The author used GPT solely for language polishing, translation, formatting, and checking the consistency of mathematical notation and references. The author reviewed, verified, and approved all AI-assisted content and takes full responsibility for the entire manuscript, including all mathematical statements and proofs.

\appendix

\section{A single spiral with a fixed angular velocity profile}\label{sec:circular}

We show that a fixed angular velocity profile gives bounded winding and vanishing dissipation under the force bound below. Fix $T>0$ and $0<a<b<\pi$. Choose $\omega\in C_c^\infty((0,\pi))$ with $\omega(r)=r-a$ on $[a,b]$. Set $w=r\omega(r)e_\varphi$ for $0<r<\pi$ and $w=0$ elsewhere on $\T^2$. Let $u^m=(v^m,\theta^m)$ be smooth solutions of \eqref{eq:NS} with the same smooth initial datum $(v_0,\theta_0)$, viscosities $\mu_m\downarrow0$ and forces $(g^m,0)$. Assume that the velocity and force are independent of $x_3$ and that

\[
 \begin{gathered}
 v^m(t)=\kappa_m'(t)w,\qquad \kappa_m\in C^\infty([0,T]),\quad \kappa_m(0)=0,\\
 \sup_m\int_0^T\|g^m(t)\|_2\,dt<\infty .
 \end{gathered}
\]

\textbf{Number of turns.} For each particle initially on the segment $\varphi=0$, $a\le r\le b$, we have $\dot r=0$ and $\dot\varphi=\kappa_m'(t)(r-a)$. Crippa, Lucà and Schulze \cite{CLS} studied mixing on the unit disk for a flow whose angular velocity is proportional to the radius. Measuring the angle continuously from $\varphi(0,r)=0$, we obtain

\begin{equation}\label{eq:single-pitch}
 \varphi(t,r)=\kappa_m(t)(r-a),\qquad
 N_m(t)=\frac{|\varphi(t,b)-\varphi(t,a)|}{2\pi}
       =\frac{|\kappa_m(t)|(b-a)}{2\pi}.
\end{equation}

At time $t$, the spiral makes $N_m(t)$ turns about the origin as $r$ increases from $a$ to $b$, including any fraction of a turn. Suppose $\kappa_m(t)\ne0$ and a fixed ray meets the spiral at least twice. From one of these points to the next, the spiral makes one full turn around the origin. The radial spacing is therefore $2\pi/|\kappa_m(t)|$.

For each $m$, let $M=\max_{0\le t\le T}\|v^m(t)\|_2$. At a time when $\|v^m(t)\|_2=M$, the horizontal energy estimate gives

\[
(M-\|v_0\|_2)(M+\|v_0\|_2)
 =M^2-\|v_0\|_2^2\le 2M\int_0^T\|g^m(s)\|_2\,ds .
\]

For $M>0$, divide by $M+\|v_0\|_2\ge M$. The resulting bound also holds when $M=0$, so

\[
|\kappa_m'(t)|\,\|w\|_2=\|v^m(t)\|_2
 \le M\le\|v_0\|_2+2\int_0^T\|g^m(s)\|_2\,ds .
\]

Since $\|w\|_2>0$, the force bound gives $\sup_m\|\kappa_m'\|_{L^\infty(0,T)}<\infty$. By \eqref{eq:single-pitch} and $v^m=\kappa_m'w$, we have

\[
 N_m(t)\le\frac{(b-a)T}{2\pi}\|\kappa_m'\|_{L^\infty(0,T)},\qquad
 \|\nabla v^m(t)\|_\infty\le\|\kappa_m'\|_{L^\infty(0,T)}\|\nabla w\|_\infty .
\]

\textbf{Dissipation estimate.} By the vertical equation, $\partial_3p^m$ is constant in $x_3$ and hence zero by periodicity of $p^m$. Thus $\theta^m$ solves \eqref{eq:scalar}. Let $\rho^m$ be the scalar transported from $\theta_0$ without diffusion. Each circle rotates through the angle $\kappa_m(t)\omega(r)$, so

\[
\rho^m(t,r,\varphi)=\theta_0\bigl(r,\varphi-\kappa_m(t)\omega(r)\bigr)
 \qquad(0<r<\pi).
\]

Outside this disk, we have $\rho^m=\theta_0$. Inside, $\partial_\varphi\rho^m$ is the same angular translate of $\partial_\varphi\theta_0$, and

\[
\partial_r\rho^m
 =(\partial_r\theta_0)\bigl(r,\varphi-\kappa_m(t)\omega(r)\bigr)
  -\kappa_m(t)\omega'(r)\,\partial_\varphi\rho^m .
\]

The angular shift preserves $r\,dr\,d\varphi$, so taking the $L^2$ norm of the gradient gives

\[
\|\nabla\rho^m(t)\|_2
 \le \bigl(1+|\kappa_m(t)|\,\|r\omega'(r)\|_\infty\bigr)\|\nabla\theta_0\|_2
 \le C_T .
\]

Here $|\kappa_m(t)|\le T\|\kappa_m'\|_\infty$, so $C_T$ is independent of $m$. Since $\rho^m(0)=\theta^m(0)$, \eqref{eq:prep} gives

\[
\mu_m\int_0^T\|\nabla\theta^m(t)\|_2^2\,dt
 \le\mu_m\int_0^T\|\nabla\rho^m(t)\|_2^2\,dt\le C_T\mu_m .
\]

Together with the bound on $\nabla v^m$ and the independence of $x_3$, this gives

\[
 \mu_m\int_0^T\|\nabla u^m(t)\|_2^2\,dt
 =\mu_m\int_0^T\bigl(\|\nabla v^m(t)\|_2^2+\|\nabla\theta^m(t)\|_2^2\bigr)\,dt
 \le C_T\mu_m\longrightarrow0 .
\]

\bibliographystyle{plain}
\begingroup
\small
\let\standardbibliography\thebibliography
\renewcommand{\thebibliography}[1]{%
  \standardbibliography{#1}\setlength{\itemsep}{1pt}\setlength{\parsep}{0pt}}
\bibliography{references}
\endgroup

\end{document}